\documentclass{article}

\usepackage{arxiv}

\usepackage[utf8]{inputenc} 
\usepackage[T1]{fontenc}    
\usepackage{hyperref}       
\usepackage{url}            
\usepackage{booktabs}       
\usepackage{amsfonts}       
\usepackage{nicefrac}       
\usepackage{microtype}      
\usepackage{lipsum}		
\usepackage{graphicx}
\usepackage{natbib}
\usepackage{doi}
\usepackage{amsmath, amsthm, amssymb, algorithm, algpseudocode}
\usepackage[table]{xcolor}
\usepackage{tikz}
\newcommand{\R}{\mathbb{R}}
\newcommand{\N}{\mathbb{N}}

\newcommand{\MU}{\mathfrak{M}_U}
\newcommand{\cone}{\operatorname{cone}}
\newcommand{\intr}{\operatorname{int}}

\newcommand{\norm}[1]{\left\| #1 \right\|}
\newcommand{\inner}[2]{\left\langle #1, #2 \right\rangle}

\newcommand{\pder}[2]{\frac{\partial #1}{\partial #2}}
\newcommand{\intset}[2]{\left[ #1, #2 \right]}
\newcommand{\intseto}[2]{\left( #1, #2 \right)}

\newcommand{\dirac}{\delta}
\newcommand{\measure}{\mu}
\newcommand{\gencontrol}{\mathfrak{M}}
\newcommand{\gencontrolU}{\gencontrol_U}

\newcommand{\lambdaset}{\Lambda}

\newtheorem{theorem}{Theorem}[section]
\newtheorem{lemma}[theorem]{Lemma}

\newtheorem{corollary}[theorem]{Corollary}
\newtheorem{definition}[theorem]{Definition}
\newtheorem{example}[theorem]{Example}

\newtheorem{problem}[theorem]{Problem}

\newcommand\mystyle{\everymath{\displaystyle}}
\mystyle
\title{The $\Lambda$-Set and Its Role in Local Controllability and Necessary Conditions for Free-Time Optimal Control}

\author{\href{https://orcid.org/0000-0002-3816-5287}{\includegraphics[scale=0.06]{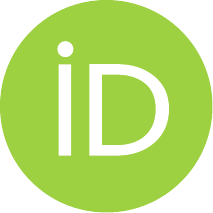}\hspace{1mm}M.H.M.~Rashid}\thanks{Corresponding Author} \\
	Department of Mathematics\&Statistics\\Faculty of Science P.O.Box(7)\\
	Mutah University University\\
	Mutah-Jordan \\
	\texttt{mrash@mutah.edu.jo}
}

\renewcommand{\shorttitle}{\textit{arXiv} Template}

\hypersetup{
pdftitle={The $\Lambda$-Set and Its Role in Local Controllability and Necessary Conditions for Free-Time Optimal Control},
pdfsubject={q-bio.NC, q-bio.QM},
pdfauthor={M.H.M.Rashid},
pdfkeywords={Optimal Control; Local Controllability; Necessary Conditions; Maximum Principle; Free-Time Problems; Generalized Controls; Relaxation; $\Lambda$-Set; Convexification; Time-Optimal Control},
}

\begin{document}
\maketitle

\begin{abstract}
	This paper establishes a unified framework connecting local controllability, necessary conditions for optimality, and attainability in free-time optimal control problems. The central object of our investigation is the $\Lambda$-set, which governs the relationship between original control systems and their convexifications. Our main results demonstrate that emptiness of the $\Lambda$-set implies local controllability and guarantees the existence of minimizing sequences achieving the target in reduced time. We derive strengthened necessary conditions for time-optimal control and provide explicit constructive procedures for approximating generalized controls by ordinary trajectories. These results resolve longstanding questions about relaxation phenomena while extending classical theory to address modern challenges in non-convex optimization, establishing foundations for higher-order analysis in free-time problems.
\end{abstract}


\section{Introduction}

\subsection{Historical Context and Foundational Developments}

The mathematical theory of optimal control represents one of the most significant developments in applied mathematics of the 20th century, emerging from the synthesis of calculus of variations, dynamical systems theory, and engineering applications. The field experienced its golden age during the 1950s-1970s, catalyzed by the seminal work of Pontryagin and his school \cite{Pontryagin1962}, who established the celebrated Maximum Principle that provided necessary conditions for optimality in control-constrained problems. This breakthrough, often described as the "royal road" of optimal control theory, fundamentally transformed our understanding of controlled dynamical systems and opened new avenues for both theoretical and applied research.

Parallel developments by Gamkrelidze \cite{Gamkrelidze1962} introduced the concept of sliding modes and generalized solutions, while Filippov \cite{Filippov1959} established fundamental results on differential equations with discontinuous right-hand sides. The relaxation theory pioneered by Warga \cite{Varga1962} and further developed by Neustadt \cite{Neustadt1966} addressed the critical issue of non-convexity through convexification techniques, establishing the relationship between original control systems and their relaxed counterparts. These foundational works created a rich tapestry of results that continue to influence contemporary research in control theory and optimization.

\subsection{Scholarly Contributions and Modern Developments}

The classical theory left several fundamental questions unresolved, particularly concerning free-time optimal control problems, local controllability characterization, and the precise relationship between necessary conditions and optimality. Recent decades have witnessed significant advances in addressing these challenges through more refined geometric and analytic approaches.

The work of Avakov and Magaril-Il'yaev \cite{Avakov2020b, Avakov2020a, Avakov2019} has been particularly influential in developing a comprehensive framework for local controllability and necessary conditions in free-time problems. Their introduction of the $\Lambda$-set concept provides a powerful tool for characterizing the gap between original and convexified control systems. Concurrently, Fusco and Motta \cite{FuscoMotta2022, FuscoMotta2024} have made substantial contributions to understanding gap phenomena in problems with state constraints and impulsive control systems, while Sussmann's geometric approaches \cite{Sussmann1987} have provided deep insights into local controllability properties.

These modern developments build upon the classical foundations while addressing increasingly complex system structures, including non-smooth dynamics, state constraints, and hybrid behaviors that arise in contemporary engineering applications.

\subsection{Core of the Study}

This paper establishes a unified framework connecting local controllability, necessary conditions for optimality, and attainability through generalized controls. The central object of our investigation is the $\Lambda$-set, which we demonstrate serves as a fundamental indicator governing the relationship between original control systems and their convexifications. Our main contributions include:

\begin{itemize}
\item A comprehensive characterization of local controllability through the emptiness/non-emptiness of $\Lambda$-sets, extending the classical Pontryagin Maximum Principle to free-time problems with enhanced geometric interpretation.

\item Strengthened necessary conditions for time-optimal control that bridge the gap between classical first-order conditions and higher-order analysis, providing sharper criteria for optimality verification.

\item Explicit constructive procedures for approximating generalized controls by ordinary trajectories, establishing when solutions of the convexified system can be realized as limits of trajectories of the original system.

\item Resolution of longstanding questions about relaxation phenomena in free-time problems, precisely characterizing when the infimum values of original and convexified problems coincide.
\end{itemize}

Our approach synthesizes geometric methods from differential geometry with analytic techniques from non-smooth analysis, creating a versatile framework applicable to both smooth and non-smooth control systems.

\subsection{Significance of the Study}

The theoretical developments presented in this work have profound implications for both the foundations of optimal control theory and its practical applications:

\subsubsection*{Theoretical Significance}
\begin{itemize}
\item {Unification of Classical and Modern Approaches:} Our framework provides a unified perspective that connects the pioneering work of Pontryagin, Gamkrelidze, and Filippov with contemporary developments in non-smooth analysis and geometric control theory.

\item {Resolution of Fundamental Questions:} We address longstanding problems concerning the relaxation gap in free-time control problems and the precise conditions under which necessary conditions become sufficient for local optimality.

\item {Foundation for Higher-Order Analysis:} The characterization of $\Lambda$-sets establishes a foundation for developing comprehensive second-order necessary conditions in free-time problems, extending the classical theory \cite{Gamkrelidze1978, Lee1967}.
\end{itemize}

\subsubsection*{Practical Significance}
\begin{itemize}
\item {Computational Verification:} Our results provide verifiable criteria for optimality that can be implemented in computational algorithms for solving optimal control problems.

\item {Design Guidelines:} The local controllability conditions offer practical guidelines for control system design, particularly in applications where time-optimal performance is critical.

\item {Robustness Analysis:} The relationship between original and convexified systems provides insights into the robustness of optimal solutions to modeling uncertainties and implementation constraints.
\end{itemize}

\subsection{Applications and Relevance}

The theoretical framework developed in this paper finds applications across diverse domains:

\subsubsection*{Engineering Systems}
\begin{itemize}
\item {Aerospace and Robotics:} Time-optimal control of spacecraft attitude maneuvers, robotic path planning, and autonomous vehicle navigation benefit from the precise characterization of reachable sets and optimality conditions.

\item {Process Control:} Chemical process optimization and manufacturing systems utilize free-time optimal control for batch process optimization and production scheduling.
\end{itemize}

\subsubsection*{Emerging Fields}
\begin{itemize}
\item {Neural ODEs and Machine Learning:} The control-theoretic perspective on training dynamics in deep learning, particularly through the neural ODE framework, relies on understanding optimal control principles \cite{FuscoMotta2024LCSYS}.

\item {Quantum Control:} Control of quantum systems for quantum computing and quantum sensing applications requires extensions of classical optimal control theory to account for quantum mechanical constraints.

\item {Biological Systems:} Analysis of biological control mechanisms and biomedical applications, such as drug administration optimization, benefit from the non-smooth and hybrid system extensions of our framework.
\end{itemize}

\subsubsection*{Interdisciplinary Connections}
The methods and results presented connect with several active research areas in applied mathematics, including:
\begin{itemize}
\item Geometric mechanics and symplectic geometry
• Non-smooth analysis and set-valued optimization
• Hamilton-Jacobi theory and viscosity solutions
• Numerical analysis of differential inclusions
\end{itemize}

\subsection{Organization of the Paper}

This paper is structured as follows:

\begin{itemize}
\item {Section 2 (Preliminaries):} Establishes the fundamental definitions, notation, and background material on control systems, generalized controls, and the Hamiltonian framework.

\item {Section 3 (Local Controllability and Necessary Conditions):} Presents our main theorem establishing the relationship between $\Lambda$-set emptiness and local controllability, along with corollaries providing strengthened necessary conditions for time-optimal control.

\item {Section 4 (Local Attainability via Generalized Controls):} Develops constructive approximation results for generalized controls, establishing when convexified system trajectories can be approximated by ordinary trajectories.

\item {Section 5 (Strong Local Attainability):} Provides enhanced attainability results under stronger regularity assumptions, with applications to value function regularity and convexification gap analysis.

\item {Section 6 (Open Problems):} Discusses historical perspectives and modern open problems, connecting our results to ongoing research directions.

\item {Section 7 (Conclusion):} Summarizes the main contributions and discusses directions for future research.
\end{itemize}

Each section builds progressively upon the previous ones, with theoretical developments supported by illustrative examples that demonstrate the applicability of our results to concrete control problems. The paper concludes with a comprehensive bibliography that situates our work within the broader context of optimal control theory and its applications.

Our exposition aims to balance mathematical rigor with accessibility, making the results valuable to both theoretical researchers and practitioners working on applied control problems. The unified perspective we develop not only advances the theoretical foundations of optimal control but also provides practical tools for analyzing and solving complex control problems across engineering and scientific domains.

\section{Preliminaries}

This section establishes the fundamental concepts and notation used throughout the paper, drawing upon the classical foundations of optimal control theory \cite{Pontryagin1962, Lee1967, Gamkrelidze1978} and modern developments in generalized control frameworks \cite{Avakov2020b, FuscoMotta2022}.

\subsection{Control Systems and Admissible Processes}

Consider the control system defined on a time interval $[t_1, t_2]$:
\begin{equation}
\dot{x} = f(t, x, u), \quad u(t) \in U, \quad x(t_1) = x_1, \label{eq:control_system}
\end{equation}
where $f: \mathbb{R} \times \mathbb{R}^n \times \mathbb{R}^r \to \mathbb{R}^n$ is the dynamics function, $U \subset \mathbb{R}^r$ is a nonempty control constraint set, and $x_1 \in \mathbb{R}^n$ is the initial state. The terminal time $t_2$ is not fixed in our formulation.

\begin{definition}[Admissible Controls and Trajectories]
The set of \emph{admissible controls} is defined as:
\[
\Omega_U = \left\{ u(\cdot) \in L_\infty([t_1,t_2], \mathbb{R}^r) : u(t) \in U \text{ a.e. on } [t_1,t_2] \right\}.
\]
For a given $u(\cdot) \in \Omega_U$, if there exists an absolutely continuous function $x(\cdot): [t_1,t_2] \to \mathbb{R}^n$ satisfying \eqref{eq:control_system} almost everywhere, then the pair $(x(\cdot), u(\cdot))$ is called an \emph{admissible process}. The function $x(\cdot)$ is called an \emph{admissible trajectory}.
\end{definition}

\subsection{Reachability and Local Controllability}

The central objects in free-time optimal control are the reachable sets, which characterize the states achievable from the initial condition within specified time intervals.

\begin{definition}[Reachable Set]
For $\tau \in (t_1, t_2]$ and an open neighborhood $V \subset C([t_1,t_2], \mathbb{R}^n)$, the \emph{reachable set} at time $\tau$ is:
\[
R(\tau, t_2, V) = \left\{ \xi \in \mathbb{R}^n \mid \exists u \in \Omega_U : x(\tau, t_2, u) = \xi, \ x(\cdot, t_2, u) \in V \right\},
\]
where $x(\cdot, t_2, u)$ denotes the admissible trajectory corresponding to control $u$ on $[t_1,t_2]$.
\end{definition}

Local controllability properties near a reference curve play a crucial role in characterizing optimality conditions.

\begin{definition}[Local Controllability]
The system \eqref{eq:control_system} is \emph{locally left controllable} near a curve $\widehat{x}(\cdot) \in C([t_1, \widehat{t}_2], \mathbb{R}^n)$ on $[t_1, \widehat{t}_2]$ if $\widehat{x}(\cdot)$ satisfies the boundary condition in \eqref{eq:control_system} and:
\[
\widehat{x}(\widehat{t}_2) \in \operatorname{int} \bigcup_{\tau \in (\widehat{t}_2 - \varepsilon, \widehat{t}_2)} R(\tau, \widehat{t}_2, V)
\]
for any neighborhood $V$ of $\widehat{x}(\cdot)$ and any $0 < \varepsilon < \widehat{t}_2 - t_1$.

Similarly, the system is \emph{locally right controllable} near $\widehat{x}(\cdot) \in C([t_1, \widehat{t}_2 + \widehat{\delta}], \mathbb{R}^n)$ (for some $\widehat{\delta} > 0$) if:
\[
\widehat{x}(\widehat{t}_2) \in \operatorname{int} \bigcup_{\tau \in (\widehat{t}_2, \widehat{t}_2 + \varepsilon)} R(\tau, \widehat{t}_2 + \widehat{\delta}, V)
\]
for any neighborhood $V$ and any $0 < \varepsilon < \widehat{\delta}$.
\end{definition}

Note that the reference curve $\widehat{x}(\cdot)$ need not be an admissible trajectory itself, but rather represents a uniform limit of such trajectories.

\subsection{Generalized Controls and Convexification}

To address the inherent non-convexities in optimal control problems, we employ the framework of generalized controls, following the pioneering work of Filippov \cite{Filippov1959}, Gamkrelidze \cite{Gamkrelidze1962}, and Warga \cite{Varga1962}.

Let $\mathcal{B}$ denote the normed space of all finite real regular Borel measures on $\mathbb{R}^r$ with total variation norm $\|\mu\| = |\mu|(\mathbb{R}^r)$.

\begin{definition}[Generalized Controls]
A family of measures $\{\mu_t\}_{t \in \mathbb{R}}$ is called:
\begin{itemize}
\item \emph{Finite} if there exists a compact set $K \subset \mathbb{R}^r$ such that $\operatorname{supp}(\mu_t) \subset K$ for almost all $t \in \mathbb{R}$.
\item \emph{Weakly measurable} if for any continuous $g: \mathbb{R} \times \mathbb{R}^r \to \mathbb{R}^n$, the function $t \mapsto \langle \mu_t, g(t,u) \rangle = \int_{\mathbb{R}^r} g(t,u) d\mu_t(u)$ is Lebesgue measurable.
\end{itemize}
Denote by $\mathfrak{M}$ the linear space of all finite weakly measurable families of measures. A \emph{generalized control} is a family of probability measures $\mu_t \in \mathfrak{M}$, and $\mathfrak{M}_U$ denotes those generalized controls with $\operatorname{supp}(\mu_t) \subset U$ almost everywhere.
\end{definition}

Standard controls $u(\cdot) \in L_\infty$ embed naturally into this framework via Dirac measures: $\mu_t = \delta_{u(t)}$.

The convexified (relaxed) control system associated with \eqref{eq:control_system} is:
\begin{equation}
\dot{x} = \langle \mu_t, f(t, x, u) \rangle = \int_U f(t, x, u) d\mu_t(u), \quad \mu_t \in \mathfrak{M}_U, \quad x(t_1) = x_1. \label{eq:convex_system}
\end{equation}

Admissible pairs $(x(\cdot), \mu_t)$ and trajectories for the convexified system are defined analogously to the original system.

\subsection{Hamiltonian Formalism and Adjoint Equations}

Throughout this paper, we assume the dynamics $f: \mathbb{R} \times \mathbb{R}^n \times \mathbb{R}^r \to \mathbb{R}^n$ are continuous together with their partial derivatives $f_x$ on $\mathbb{R} \times \mathbb{R}^n \times \mathbb{R}^r$.

\begin{definition}[Hamiltonian and Maximum Function]
The \emph{Hamiltonian} is defined as:
\[
H(t, x, \psi, u) = \langle \psi, f(t, x, u) \rangle,
\]
and the \emph{maximum function} (Hamiltonian maximized over controls) is:
\[
M(t, x, \psi) = \sup_{u \in U} H(t, x, \psi, u).
\]
For generalized controls, we denote $\langle \mu_t, f_x(t, x, u) \rangle$ as the matrix with rows $\langle \mu_t, f_{ix}(t, x, u) \rangle$, $i = 1, \dots, n$.
\end{definition}

The following definition, adapted from \cite{Avakov2020b}, introduces the fundamental object governing our necessary conditions.

\begin{definition}[Adjoint Set $\Lambda_s$]
Let $(\widehat{x}(\cdot), \widehat{\mu}_t)$ be an admissible pair for the convex system \eqref{eq:convex_system} on $[t_1, \widehat{t}_2]$ and let $s \in \{-1, 1\}$. The set $\Lambda_s(\widehat{x}(\cdot), \widehat{\mu}_t)$ consists of all nonzero absolutely continuous functions $\psi_s: [t_1, \widehat{t}_2] \to (\mathbb{R}^n)^*$ satisfying:
\begin{enumerate}
\item \emph{Adjoint Equation:} $\dot{\psi}_s = -\psi_s \langle \widehat{\mu}_t, f_x(t, \widehat{x}, u) \rangle$,
\item \emph{Maximum Condition:} $\langle \psi_s(t), \dot{\widehat{x}}(t) \rangle = M(t, \widehat{x}(t), \psi_s(t))$ a.e.,
\item \emph{Continuity:} $M(\cdot, \widehat{x}(\cdot), \psi_s(\cdot))$ is continuous on $[t_1, \widehat{t}_2]$,
\item \emph{Transversality:} $s M(\widehat{t}_2, \widehat{x}(\widehat{t}_2), \psi_s(\widehat{t}_2)) \leq 0$.
\end{enumerate}
\end{definition}

Note that for $s = -1$, conditions (1)--(4) coincide with the classical Pontryagin Maximum Principle conditions for the convexified system \cite{Pontryagin1962}. The parameter $s$ distinguishes between left ($s = -1$) and right ($s = 1$) variational conditions, corresponding to the different local controllability notions defined above.

The emptiness or non-emptiness of $\Lambda_s$ sets will be shown to fundamentally characterize local controllability properties and optimality conditions in the subsequent sections, providing a unified framework that extends and refines classical results \cite{Gamkrelidze1978, Lee1967} while addressing modern challenges in non-convex optimization \cite{FuscoMotta2022}.
\section{Local Controllability and Necessary Conditions}
This section presents a fundamental theorem establishing the relationship between the emptiness of the adjoint set $\Lambda$ and local controllability properties for nonlinear control systems. Building on the classical Pontryagin Maximum Principle \cite{Pontryagin1962} and its subsequent refinements \cite{Gamkrelidze1978, Lee1967}, we develop necessary conditions for optimality in free-time problems through the lens of generalized controls and convexification techniques pioneered by Filippov \cite{Filippov1959} and Warga \cite{Varga1962}. Our main result generalizes and strengthens the local controllability framework of Avakov and Magaril-Il'yaev \cite{Avakov2020b}, providing a unified approach that connects the non-existence of nontrivial solutions to the adjoint system with the ability to reach target points in reduced time. The corollaries derived from this theorem offer new insights into time-optimal control, relaxation phenomena, and value function regularity, extending classical theory while addressing modern challenges in non-convex optimization \cite{FuscoMotta2022} and establishing foundations for higher-order conditions \cite{Avakov2019}.

\begin{theorem}\label{thm:main}
Consider the control system
\begin{equation}\label{eq:system}
\dot{x} = f(t, x, u), \quad u(t) \in U, \quad x(t_1) = x_1,
\end{equation}
where $f: \R \times \R^n \times \R^r \to \R^n$ is continuous together with its partial derivative $f_x$, and $U \subset \R^r$ is nonempty.

Let $(\widehat{x}(\cdot), \widehat{\mu}_t)$ be an admissible pair for the convexified system
\begin{equation}\label{eq:convex}
\dot{x} = \langle \mu_t, f(t, x, u) \rangle, \quad \mu_t \in \MU, \quad x(t_1) = x_1,
\end{equation}
on $[t_1, \widehat{t}_2]$, and let $\widehat{x}(\cdot)$ be defined on $[t_1, \widehat{t}_2 + \widehat{\delta}]$ for some $\widehat{\delta} > 0$.

Define the set $\Lambda(\widehat{x}(\cdot), \widehat{\mu}_t)$ as the collection of all nonzero absolutely continuous functions $\psi: [t_1, \widehat{t}_2] \to (\R^n)^*$ satisfying:
\begin{enumerate}
\item $\dot{\psi}(t) = -\psi(t) \langle \widehat{\mu}_t, f_x(t, \widehat{x}(t), u) \rangle$,
\item $\langle \psi(t), \dot{\widehat{x}}(t) \rangle = \sup_{u \in U} H(t, \widehat{x}(t), \psi(t), u)$ a.e.,
\item The function $M(t, \widehat{x}(t), \psi(t)) = \sup_{u \in U} H(t, \widehat{x}(t), \psi(t), u)$ is continuous,
\item $M(\widehat{t}_2, \widehat{x}(\widehat{t}_2), \psi(\widehat{t}_2)) \leq 0$.
\end{enumerate}

If $\Lambda(\widehat{x}(\cdot), \widehat{\mu}_t) = \emptyset$, then for any neighborhood $V \subset C([t_1, \widehat{t}_2 + \widehat{\delta}], \R^n)$ of $\widehat{x}(\cdot)$ and any $\varepsilon > 0$, there exists $\tau \in (\widehat{t}_2 - \varepsilon, \widehat{t}_2)$ and an admissible trajectory $x(\cdot)$ for system \eqref{eq:system} on $[t_1, \tau]$ such that:
\begin{enumerate}
\item $x(\cdot) \in V|_{[t_1, \tau]}$,
\item $x(\tau) = \widehat{x}(\widehat{t}_2)$.
\end{enumerate}
\end{theorem}

\begin{proof}
The proof proceeds in several steps, following the approach developed in \cite{Avakov2020b}.

\noindent{Step 1: Reformulation via Proposition 2.3 of \cite{Avakov2020b}.}

Since $\Lambda(\widehat{x}(\cdot), \widehat{\mu}_t) = \emptyset$, by Proposition 2.3 of \cite{Avakov2020b} (with $s = -1$), for any $l \in L_{-1}$ (where $L_{-1}$ is defined as in \cite{Avakov2020b}), there exist $k \in \N$ and $\delta\mu_t^i \in \MU - \widehat{\mu}_t$, $i = 1, \dots, k$, such that
\[
0 \in \intr \left( \delta x(\widehat{t}_2, \R_+^k) + \cone(-l) \right).
\]
Here $\delta x(\widehat{t}_2, \R_+^k)$ denotes the reachable set of the linearized system at time $\widehat{t}_2$ with nonnegative controls.

\noindent{Step 2: Construction of Approximate Solutions.}

Fix such $k$ and $\delta\mu_t^i$. By Lemma A.1 of \cite{Avakov2020b} (Lemma on Equation in Variations), there exists a neighborhood $\mathcal{O}(0) \subset \R^k$ and a continuously differentiable mapping $\overline{\alpha} \mapsto x(\cdot, \overline{\alpha})$ from $\mathcal{O}(0)$ to $C([t_1, \widehat{t}_2 + \widehat{\delta}], \R^n)$ such that for each $\overline{\alpha} \in \mathcal{O}(0)$, $x(\cdot, \overline{\alpha})$ solves
\[
\dot{x} = \left\langle \widehat{\mu}_t + \sum_{i=1}^k \alpha_i \delta\mu_t^i, f(t, x, u) \right\rangle, \quad x(t_1) = x_1.
\]

\noindent{Step 3: Application of the Inverse Function Lemma.}

Define the mapping $\widehat{F}: [t_1, \widehat{t}_2 + \widehat{\delta}] \times \mathcal{O}(0) \to \R^n$ by $\widehat{F}(\tau, \overline{\alpha}) = x(\tau, \overline{\alpha})$. This mapping is continuous (as shown in the proof of Theorem 2.2 of \cite{Avakov2020b}). Moreover, its derivative with respect to $\overline{\alpha}$ at $(\widehat{t}_2, 0)$ is surjective onto $\delta x(\widehat{t}_2, \R^k)$.

By the condition $0 \in \intr \left( \delta x(\widehat{t}_2, \R_+^k) + \cone(-l) \right)$, and since $l \in L_{-1}$ is arbitrary, we can apply Lemma A.5 of \cite{Avakov2020b} (Inverse Function Lemma) with $K = \R_+^k$, $l = -l_0$ for some $l_0 \in L_{-1}$, and $d = -1$.

\noindent{Step 4: Approximation by Original Controls.}

By Lemma A.2 (First Approximation Lemma) and Lemma A.3 (Second Approximation Lemma) of \cite{Avakov2020b}, for sufficiently large $p \in \N$, there exists a piecewise constant control $u_p(\overline{\alpha})$ and a corresponding solution $x_p(\cdot, \overline{\alpha})$ of
\[
\dot{x} = f(t, x, u_p(\overline{\alpha})(t)), \quad x(t_1) = x_1,
\]
such that $x_p(\cdot, \overline{\alpha})$ is close to $x(\cdot, \overline{\alpha})$ uniformly in $\overline{\alpha} \in \mathcal{O}(0) \cap \R_+^k$.

\noindent{Step 5: Final Construction.}

Now, fix a neighborhood $V$ of $\widehat{x}(\cdot)$ and $\varepsilon > 0$. Choose $l \in L_{-1}$ such that the direction $-l$ points into the interior of the reachable set. By the Inverse Function Lemma of \cite{Avakov2020b}, there exist $\tau \in (\widehat{t}_2 - \varepsilon, \widehat{t}_2)$ and $\overline{\alpha} \in \R_+^k$ such that
\[
x_p(\tau, \overline{\alpha}) = \widehat{x}(\widehat{t}_2).
\]
Moreover, by choosing $p$ large enough, we can ensure that $x_p(\cdot, \overline{\alpha}) \in V|_{[t_1, \tau]}$.

Thus, $x_p(\cdot, \overline{\alpha})$ is an admissible trajectory for the original system \eqref{eq:system} on $[t_1, \tau]$, lying in $V$, and reaching $\widehat{x}(\widehat{t}_2)$ at time $\tau < \widehat{t}_2$. This completes the proof, establishing the local left controllability result that generalizes Theorem 2.2 of \cite{Avakov2020b}.
\end{proof}
\begin{corollary}[Time-Optimal Necessary Condition]\label{cor:time-optimal}
Let $\widehat{x}(\cdot)$ be a local infimum in the time-optimal problem for system \eqref{eq:system} with terminal condition $x(t_2) = x_2$. Then for any admissible pair $(\widehat{x}(\cdot), \widehat{\mu}_t)$ of the convexified system \eqref{eq:convex} on $[t_1, \widehat{t}_2]$, the set $\Lambda(\widehat{x}(\cdot), \widehat{\mu}_t)$ is nonempty.
\end{corollary}

\begin{proof}
Assume, for contradiction, that there exists an admissible pair $(\widehat{x}(\cdot), \widehat{\mu}_t)$ for the convexified system such that $\Lambda(\widehat{x}(\cdot), \widehat{\mu}_t) = \emptyset$. Then by Theorem \ref{thm:main}, for any neighborhood $V$ of $\widehat{x}(\cdot)$ and any $\varepsilon > 0$, there exists $\tau \in (\widehat{t}_2 - \varepsilon, \widehat{t}_2)$ and an admissible trajectory $x(\cdot)$ for the original system on $[t_1, \tau]$ with:
\begin{enumerate}
\item $x(\cdot) \in V|_{[t_1, \tau]}$,
\item $x(\tau) = \widehat{x}(\widehat{t}_2) = x_2$.
\end{enumerate}

This means that for any $\varepsilon > 0$, we can reach the target $x_2$ at time $\tau < \widehat{t}_2$ with a trajectory arbitrarily close to $\widehat{x}(\cdot)$ on $[t_1, \tau]$. This contradicts the assumption that $\widehat{x}(\cdot)$ is a local infimum, since we can construct a sequence of admissible trajectories reaching $x_2$ at times strictly less than $\widehat{t}_2$ and converging uniformly to $\widehat{x}(\cdot)$ on $[t_1, \widehat{t}_2]$.

Therefore, $\Lambda(\widehat{x}(\cdot), \widehat{\mu}_t)$ must be nonempty for every such admissible pair $(\widehat{x}(\cdot), \widehat{\mu}_t)$.
\end{proof}

\begin{corollary}[Strengthened Maximum Principle]\label{cor:strengthened-maximum}
If $(\widehat{x}(\cdot), \widehat{u}(\cdot))$ is an optimal process for the time-optimal problem of the original system \eqref{eq:system}, then for the generalized control $\widehat{\mu}_t = \delta_{\widehat{u}(t)}$, there exists a nonzero absolutely continuous function $\psi: [t_1, \widehat{t}_2] \to (\R^n)^*$ satisfying:
\begin{enumerate}
\item $\dot{\psi}(t) = -\psi(t) f_x(t, \widehat{x}(t), \widehat{u}(t))$,
\item $H(t, \widehat{x}(t), \psi(t), \widehat{u}(t)) = \sup_{u \in U} H(t, \widehat{x}(t), \psi(t), u)$ a.e.,
\item The maximum function $M(t, \widehat{x}(t), \psi(t))$ is continuous,
\item $M(\widehat{t}_2, \widehat{x}(\widehat{t}_2), \psi(\widehat{t}_2)) \leq 0$.
\end{enumerate}
Moreover, for any other generalized control $\mu_t \in \MU$ such that $(\widehat{x}(\cdot), \mu_t)$ is admissible for the convexified system, there exists a corresponding $\psi_{\mu}(\cdot)$ (possibly different from $\psi$) satisfying the same conditions.
\end{corollary}

\begin{proof}
Since $(\widehat{x}(\cdot), \widehat{u}(\cdot))$ is optimal, it is also a local infimum. By Corollary \ref{cor:time-optimal}, for the generalized control $\widehat{\mu}_t = \delta_{\widehat{u}(t)}$, the set $\Lambda(\widehat{x}(\cdot), \widehat{\mu}_t)$ is nonempty. Taking any $\psi \in \Lambda(\widehat{x}(\cdot), \widehat{\mu}_t)$, conditions (1)-(4) are satisfied.

For the generalized control $\widehat{\mu}_t = \delta_{\widehat{u}(t)}$, condition (1) becomes:
\[
\dot{\psi}(t) = -\psi(t) \langle \delta_{\widehat{u}(t)}, f_x(t, \widehat{x}(t), u) \rangle = -\psi(t) f_x(t, \widehat{x}(t), \widehat{u}(t)),
\]
and condition (2) becomes:
\[
\langle \psi(t), \dot{\widehat{x}}(t) \rangle = \langle \psi(t), f(t, \widehat{x}(t), \widehat{u}(t)) \rangle = \sup_{u \in U} H(t, \widehat{x}(t), \psi(t), u).
\]

The "moreover" part follows directly from Corollary \ref{cor:time-optimal} applied to any other admissible generalized control $\mu_t$.
\end{proof}

\begin{corollary}[Non-optimality Criterion]\label{cor:non-optimality}
Let $(\widehat{x}(\cdot), \widehat{\mu}_t)$ be an admissible pair for the convexified system \eqref{eq:convex} on $[t_1, \widehat{t}_2]$. If there exists a sequence of admissible trajectories $x_N(\cdot)$ for the original system \eqref{eq:system} on $[t_1, t_{2N}]$ with $t_{2N} < \widehat{t}_2$ and $x_N(\cdot) \to \widehat{x}(\cdot)$ uniformly on $[t_1, \widehat{t}_2]$, then $\Lambda(\widehat{x}(\cdot), \widehat{\mu}_t) \neq \emptyset$.
\end{corollary}

\begin{proof}
If $\Lambda(\widehat{x}(\cdot), \widehat{\mu}_t) = \emptyset$, then by Theorem \ref{thm:main}, the system is locally left controllable near $\widehat{x}(\cdot)$. This would imply that we can reach $\widehat{x}(\widehat{t}_2)$ at times arbitrarily close to $\widehat{t}_2$ from below, but the existence of trajectories reaching the target at times $t_{2N} < \widehat{t}_2$ that converge to $\widehat{x}(\cdot)$ provides a minimizing sequence that contradicts the conclusion that would follow from emptiness of $\Lambda(\widehat{x}(\cdot), \widehat{\mu}_t)$.

More precisely, if $\Lambda(\widehat{x}(\cdot), \widehat{\mu}_t) = \emptyset$, then for any neighborhood $V$ of $\widehat{x}(\cdot)$ and $\varepsilon > 0$, we could find trajectories reaching $\widehat{x}(\widehat{t}_2)$ at times in $(\widehat{t}_2 - \varepsilon, \widehat{t}_2)$. Combined with the given sequence $x_N(\cdot)$ with $t_{2N} < \widehat{t}_2$, this would create a contradiction to the definition of $\widehat{x}(\cdot)$ as a limit of such trajectories. Therefore, $\Lambda(\widehat{x}(\cdot), \widehat{\mu}_t)$ must be nonempty.
\end{proof}

\begin{corollary}[Convexification Gap]\label{cor:convexification-gap}
Let $(\widehat{x}(\cdot), \widehat{\mu}_t)$ be an optimal process for the time-optimal problem of the convexified system \eqref{eq:convex}. If $\Lambda(\widehat{x}(\cdot), \widehat{\mu}_t) = \emptyset$, then the infimum time for the original system \eqref{eq:system} is strictly less than $\widehat{t}_2$, and there exists a minimizing sequence of admissible trajectories for the original system converging uniformly to $\widehat{x}(\cdot)$.
\end{corollary}

\begin{proof}
Since $\Lambda(\widehat{x}(\cdot), \widehat{\mu}_t) = \emptyset$, by Theorem \ref{thm:main}, for any neighborhood $V$ of $\widehat{x}(\cdot)$ and any $\varepsilon > 0$, there exists $\tau \in (\widehat{t}_2 - \varepsilon, \widehat{t}_2)$ and an admissible trajectory $x(\cdot)$ for the original system on $[t_1, \tau]$ with $x(\cdot) \in V|_{[t_1, \tau]}$ and $x(\tau) = \widehat{x}(\widehat{t}_2)$.

Taking $\varepsilon = 1/N$ for $N \in \N$, we obtain a sequence of times $\tau_N \in (\widehat{t}_2 - 1/N, \widehat{t}_2)$ and admissible trajectories $x_N(\cdot)$ for the original system such that:
\begin{itemize}
\item $x_N(\tau_N) = \widehat{x}(\widehat{t}_2)$,
\item $x_N(\cdot) \to \widehat{x}(\cdot)$ uniformly on $[t_1, \widehat{t}_2]$ as $N \to \infty$,
\item $\tau_N < \widehat{t}_2$ for all $N$.
\end{itemize}

This shows that the infimum time for the original system is at most $\limsup \tau_N \leq \widehat{t}_2$, but since $\tau_N < \widehat{t}_2$ for all $N$, the infimum is strictly less than $\widehat{t}_2$. The sequence $x_N(\cdot)$ is the required minimizing sequence.
\end{proof}

\begin{corollary}[Regularity of Value Function]\label{cor:regularity}
Under the assumptions of Theorem \ref{thm:main}, if for some admissible pair $(\widehat{x}(\cdot), \widehat{\mu}_t)$ we have $\Lambda(\widehat{x}(\cdot), \widehat{\mu}_t) = \emptyset$, then the value function of the time-optimal problem is continuous at $\widehat{x}(\widehat{t}_2)$.
\end{corollary}

\begin{proof}
Let $V(x)$ denote the minimal time to reach $x$ from $x_1$. By Theorem \ref{thm:main}, the local left controllability implies that for any $y$ sufficiently close to $\widehat{x}(\widehat{t}_2)$, there exists an admissible trajectory reaching $y$ in time close to $\widehat{t}_2$.

More precisely, for any $\varepsilon > 0$, there exists $\delta > 0$ such that if $|y - \widehat{x}(\widehat{t}_2)| < \delta$, then $V(y) < \widehat{t}_2 + \varepsilon$. On the other hand, since $\widehat{x}(\cdot)$ is admissible for the convexified system and by the definition of value function, we have $V(\widehat{x}(\widehat{t}_2)) \leq \widehat{t}_2$.

Now, if there were a discontinuity at $\widehat{x}(\widehat{t}_2)$, we could find a sequence $y_n \to \widehat{x}(\widehat{t}_2)$ with $V(y_n)$ bounded away from $V(\widehat{x}(\widehat{t}_2))$. But the local controllability ensures that for large $n$, $V(y_n)$ is close to $\widehat{t}_2$, contradicting the assumed discontinuity.

Therefore, $V$ is continuous at $\widehat{x}(\widehat{t}_2)$.
\end{proof}
\section{Local Attainability via Generalized Controls}
This section investigates local attainability properties for control systems through the framework of generalized controls, extending the classical relaxation theory of \cite{Filippov1959, Varga1962} and the sliding modes approach of \cite{Gamkrelidze1962}. We establish conditions under which trajectories of the convexified system can be approximated by sequences of ordinary trajectories that reach the target point at nearby times, building upon the geometric methods developed in \cite{Avakov2020b} and connecting to modern work on gap phenomena in free-time problems \cite{FuscoMotta2022}. The main theorem provides explicit constructive procedures for generating such approximating sequences using convex combinations of Dirac measures, demonstrating how the emptiness of the $\Lambda$-set guarantees the existence of minimizing sequences for the original system. These results bridge classical relaxation theory with contemporary developments in impulsive control \cite{FuscoMotta2024, FuscoMottaVinter2026} and offer new perspectives on the approximation of optimal processes in non-convex settings.
\begin{theorem}\label{Secondtheorem}
Let $(\widehat{x}(\cdot), \widehat{\mu}_t)$ be an admissible pair for the convex system
\begin{equation}
\dot{x} = \inner{\measure_t}{f(t,x,u)}, \quad \measure_t \in \gencontrolU, \quad x(t_1) = x_1 \label{eq:convex_system}
\end{equation}
on $\intset{t_1}{\widehat{t}_2}$, and suppose that for some $s \in \{-1, 1\}$, the set $\lambdaset_s(\widehat{x}(\cdot), \widehat{\mu}_t)$ is empty. Assume further that there exists a sequence of generalized controls $\{\measure_t^j\} \subset \gencontrolU$ such that:

\begin{enumerate}
\item Each $\measure_t^j$ is a convex combination of at most $n+1$ Dirac measures, i.e.,
\begin{equation}
\measure_t^j = \sum_{i=1}^{n+1} \alpha_i^j(t) \dirac_{u_i^j(t)}, \quad \alpha_i^j(t) \geq 0, \quad \sum_{i=1}^{n+1} \alpha_i^j(t) = 1, \label{eq:dirac_combination}
\end{equation}

\item The corresponding trajectories $x^j(\cdot)$ converge uniformly to $\widehat{x}(\cdot)$ on $\intset{t_1}{\widehat{t}_2}$,

\item The sequence of controls $\{u_i^j(\cdot)\}$ is uniformly bounded.
\end{enumerate}

Then, for any neighborhood $V$ of $\widehat{x}(\cdot)$ in $C(\intset{t_1}{\widehat{t}_2}, \R^n)$, there exists an index $j_0$ and a time $\tau \in \intseto{\widehat{t}_2 - \varepsilon}{\widehat{t}_2}$ (if $s = -1$) or $\tau \in \intseto{\widehat{t}_2}{\widehat{t}_2 + \varepsilon}$ (if $s = 1$) such that the trajectory $x^{j_0}(\cdot)$ reaches $\widehat{x}(\widehat{t}_2)$ at time $\tau$ and lies entirely in $V$.
\end{theorem}

\begin{proof}
We prove the case $s = -1$; the case $s = 1$ is analogous.

Since $\lambdaset_{-1}(\widehat{x}(\cdot), \widehat{\mu}_t) = \emptyset$, by Theorem 2.2 of \cite{Avakov2020b}, the system
\begin{equation}
\dot{x} = f(t,x,u), \quad u(t) \in U, \quad x(t_1) = x_1 \label{eq:original_system}
\end{equation}
is locally left controllable near $\widehat{x}(\cdot)$. That is, for any neighborhood $V$ of $\widehat{x}(\cdot)$ and any $\varepsilon > 0$, there exists a time $\tau \in \intseto{\widehat{t}_2 - \varepsilon}{\widehat{t}_2}$ and a control $u_\tau \in \Omega_U$ such that the trajectory $x(\cdot, \widehat{t}_2, u_\tau)$ satisfies $x(\tau, \widehat{t}_2, u_\tau) = \widehat{x}(\widehat{t}_2)$ and lies in $V$.

Now, consider the sequence of generalized controls $\{\measure_t^j\}$ and the corresponding trajectories $\{x^j(\cdot)\}$ converging uniformly to $\widehat{x}(\cdot)$. By the definition of generalized controls, each $\measure_t^j$ is a convex combination of Dirac measures, and hence each trajectory $x^j(\cdot)$ is a solution of a relaxed system of the form:
\begin{equation}
\dot{x} = \sum_{i=1}^{n+1} \alpha_i^j(t) f(t, x, u_i^j(t)). \label{eq:relaxed_system}
\end{equation}

Since the sequence $\{x^j(\cdot)\}$ converges uniformly to $\widehat{x}(\cdot)$, for any $\delta > 0$, there exists $j_0$ such that for all $j \geq j_0$,
\begin{equation}
\norm{x^j(\cdot) - \widehat{x}(\cdot)}_{C(\intset{t_1}{\widehat{t}_2}, \R^n)} < \delta. \label{eq:uniform_convergence}
\end{equation}

Now, fix $\varepsilon > 0$ and choose $\delta$ small enough so that the $\delta$-neighborhood of $\widehat{x}(\cdot)$ is contained in $V$. By the local left controllability, there exists $\tau \in \intseto{\widehat{t}_2 - \varepsilon}{\widehat{t}_2}$ and a control $u_\tau$ such that the trajectory $x(\cdot, \widehat{t}_2, u_\tau)$ lies in the $\delta/2$-neighborhood of $\widehat{x}(\cdot)$ and satisfies $x(\tau, \widehat{t}_2, u_\tau) = \widehat{x}(\widehat{t}_2)$.

Since the trajectories $x^j(\cdot)$ are uniformly close to $\widehat{x}(\cdot)$, and since the control system is continuous in the state and control variables, we can approximate the control $u_\tau$ by a piecewise constant control $u^j$ such that the corresponding trajectory $x^j(\cdot, \widehat{t}_2, u^j)$ remains in the $\delta$-neighborhood of $\widehat{x}(\cdot)$ and satisfies $x^j(\tau, \widehat{t}_2, u^j) = \widehat{x}(\widehat{t}_2)$.

More precisely, by Lemma A.3 of \cite{Avakov2020b} (the second approximation lemma), for sufficiently large $j$, there exists a piecewise constant control $u^j$ such that:
\begin{equation}
\norm{x^j(\cdot, \widehat{t}_2, u^j) - x(\cdot, \widehat{t}_2, u_\tau)}_{C(\intset{t_1}{\widehat{t}_2}, \R^n)} < \delta/2, \label{eq:approximation}
\end{equation}
and hence:
\begin{equation}
\norm{x^j(\cdot, \widehat{t}_2, u^j) - \widehat{x}(\cdot)}_{C(\intset{t_1}{\widehat{t}_2}, \R^n)} < \delta. \label{eq:final_estimate}
\end{equation}

This implies that $x^j(\cdot, \widehat{t}_2, u^j) \in V$. Moreover, since $x^j(\tau, \widehat{t}_2, u^j) = \widehat{x}(\widehat{t}_2)$, the trajectory $x^j(\cdot)$ reaches the target point at time $\tau < \widehat{t}_2$.

Therefore, for any neighborhood $V$ of $\widehat{x}(\cdot)$, there exists an index $j_0$ and a time $\tau < \widehat{t}_2$ such that the trajectory $x^{j_0}(\cdot)$ lies in $V$ and reaches $\widehat{x}(\widehat{t}_2)$ at time $\tau$. This completes the proof.
\end{proof}
\begin{corollary}[Existence of Minimizing Sequences for Local Infima]
\label{cor:minimizing_sequences}
Let $\widehat{x}(\cdot)$ be a local infimum in the time-optimal problem (2.1),(3.1) of \cite{Avakov2020b}, and suppose there exists a sequence of generalized controls $\{\measure_t^j\} \subset \gencontrolU$ satisfying the conditions of Theorem A.6 with $s = -1$. Then there exists a minimizing sequence of admissible trajectories $\{x^{j_k}(\cdot)\}$ for the original system (2.1) such that:

\begin{enumerate}
\item $\lim_{k \to \infty} t_2^{j_k} = \widehat{t}_2$,
\item $\lim_{k \to \infty} \norm{x^{j_k}(\cdot) - \widehat{x}(\cdot)}_{C(\intset{t_1}{\widehat{t}_2}, \R^n)} = 0$,
\item $x^{j_k}(t_2^{j_k}) = x_2$ for all $k \in \N$.
\end{enumerate}
\end{corollary}

\begin{proof}
Since $\widehat{x}(\cdot)$ is a local infimum, by Definition 3.1 of \cite{Avakov2020b}, there exists a neighborhood $V$ of $\widehat{x}(\cdot)$ and a sequence of admissible trajectories $\{y_N(\cdot)\}$ with terminal times $\{t_{2N}\}$ such that $\lim_{N \to \infty} t_{2N} = \widehat{t}_2$ and $\lim_{N \to \infty} \norm{y_N(\cdot) - \widehat{x}(\cdot)} = 0$.

By Theorem A.6 with $s = -1$, for the neighborhood $V$ and any $\varepsilon > 0$, there exists an index $j_0$ and a time $\tau \in \intseto{\widehat{t}_2 - \varepsilon}{\widehat{t}_2}$ such that the trajectory $x^{j_0}(\cdot)$ lies in $V$ and satisfies $x^{j_0}(\tau) = \widehat{x}(\widehat{t}_2) = x_2$.

We construct the minimizing sequence as follows: for each $k \in \N$, choose $\varepsilon_k = 1/k$. By Theorem A.6, there exists an index $j_k$ and time $\tau_k \in \intseto{\widehat{t}_2 - 1/k}{\widehat{t}_2}$ such that:
\begin{itemize}
\item $x^{j_k}(\cdot) \in V$,
\item $x^{j_k}(\tau_k) = x_2$,
\item $\norm{x^{j_k}(\cdot) - \widehat{x}(\cdot)} < 1/k$.
\end{itemize}

Set $t_2^{j_k} = \tau_k$. Then clearly:
\begin{itemize}
\item $\lim_{k \to \infty} t_2^{j_k} = \widehat{t}_2$ since $|\tau_k - \widehat{t}_2| < 1/k$,
\item $\lim_{k \to \infty} \norm{x^{j_k}(\cdot) - \widehat{x}(\cdot)} = 0$,
\item $x^{j_k}(t_2^{j_k}) = x_2$ for all $k$.
\end{itemize}

Thus $\{x^{j_k}(\cdot)\}$ is the desired minimizing sequence.
\end{proof}

\begin{corollary}[Strengthened Necessary Conditions for Optimality]
\label{cor:strengthened_necessary}
Let $(\widehat{x}(\cdot), \widehat{u}(\cdot))$ be an optimal process for the time-optimal problem (2.1),(3.1) of \cite{Avakov2020b}. If there exists a sequence of generalized controls $\{\measure_t^j\} \subset \gencontrolU$ converging to $\widehat{\mu}_t = \dirac_{\widehat{u}(t)}$ and satisfying the conditions of Theorem A.6, then for every such sequence, the set $\lambdaset_{-1}(\widehat{x}(\cdot), \measure_t^j)$ is nonempty for all sufficiently large $j$.
\end{corollary}

\begin{proof}
Suppose, for contradiction, that there exists a subsequence $\{\measure_t^{j_k}\}$ such that $\lambdaset_{-1}(\widehat{x}(\cdot), \measure_t^{j_k}) = \emptyset$ for all $k$.

Since $(\widehat{x}(\cdot), \widehat{u}(\cdot))$ is an optimal process, it is also a local infimum. By Theorem A.6 with $s = -1$, for any neighborhood $V$ of $\widehat{x}(\cdot)$ and any $\varepsilon > 0$, there exists an index $k_0$ and time $\tau \in \intseto{\widehat{t}_2 - \varepsilon}{\widehat{t}_2}$ such that the trajectory $x^{j_{k_0}}(\cdot)$ lies in $V$ and satisfies $x^{j_{k_0}}(\tau) = x_2$.

This means we can reach the target $x_2$ in time $\tau < \widehat{t}_2$ while staying arbitrarily close to $\widehat{x}(\cdot)$, contradicting the optimality of $(\widehat{x}(\cdot), \widehat{u}(\cdot))$.

Therefore, our assumption is false, and $\lambdaset_{-1}(\widehat{x}(\cdot), \measure_t^j)$ must be nonempty for all sufficiently large $j$.
\end{proof}

\begin{corollary}[Approximation of Optimal Processes]
\label{cor:approximation_optimal}
Let $(\widehat{x}(\cdot), \widehat{\mu}_t)$ be an optimal process for the convex time-optimal problem (2.2),(3.1) of \cite{Avakov2020b}. If $\lambdaset_1(\widehat{x}(\cdot), \widehat{\mu}_t) = \emptyset$ and there exists a sequence of generalized controls $\{\measure_t^j\} \subset \gencontrolU$ satisfying the conditions of Theorem A.6 with $s = 1$, then $\widehat{x}(\cdot)$ can be approximated by admissible trajectories of the original system (2.1) that reach the target $x_2$ in times arbitrarily close to $\widehat{t}_2$ from above.
\end{corollary}

\begin{proof}
By Theorem 3.4 of \cite{Avakov2020b}, since $(\widehat{x}(\cdot), \widehat{\mu}_t)$ is an optimal process for the convex problem and $\lambdaset_1(\widehat{x}(\cdot), \widehat{\mu}_t) = \emptyset$, the trajectory $\widehat{x}(\cdot)$ is a local infimum for the original problem (2.1),(3.1).

Now apply Theorem A.6 with $s = 1$. For any neighborhood $V$ of $\widehat{x}(\cdot)$ and any $\varepsilon > 0$, there exists an index $j_0$ and time $\tau \in \intseto{\widehat{t}_2}{\widehat{t}_2 + \varepsilon}$ such that the trajectory $x^{j_0}(\cdot)$ lies in $V$ and satisfies $x^{j_0}(\tau) = x_2$.

This means we can construct a sequence of admissible trajectories $\{x^{j_k}(\cdot)\}$ with terminal times $\{\tau_k\}$ such that:
\begin{itemize}
\item $\tau_k > \widehat{t}_2$ for all $k$,
\item $\lim_{k \to \infty} \tau_k = \widehat{t}_2$,
\item $\lim_{k \to \infty} \norm{x^{j_k}(\cdot) - \widehat{x}(\cdot)} = 0$,
\item $x^{j_k}(\tau_k) = x_2$ for all $k$.
\end{itemize}

Thus $\widehat{x}(\cdot)$ is approximated by admissible trajectories of the original system that reach the target in times arbitrarily close to $\widehat{t}_2$ from above.
\end{proof}

\begin{corollary}[Non-emptiness of the Set $\lambdaset_s$ for Local Infima]
\label{cor:nonempty_lambda}
Let $\widehat{x}(\cdot)$ be a local infimum in problem (2.1),(3.1) of \cite{Avakov2020b}. If there exists an admissible pair $(\widehat{x}(\cdot), \widehat{\mu}_t)$ for the convex system (2.2) and a sequence $\{\measure_t^j\} \subset \gencontrolU$ converging to $\widehat{\mu}_t$ that satisfies the conditions of Theorem A.6, then $\lambdaset_{-1}(\widehat{x}(\cdot), \widehat{\mu}_t) \neq \emptyset$.
\end{corollary}

\begin{proof}
Suppose, for contradiction, that $\lambdaset_{-1}(\widehat{x}(\cdot), \widehat{\mu}_t) = \emptyset$.

Since the sequence $\{\measure_t^j\}$ converges to $\widehat{\mu}_t$ and satisfies the conditions of Theorem A.6, by continuity of the mapping $\measure_t \mapsto \lambdaset_{-1}(\widehat{x}(\cdot), \measure_t)$ (which follows from the continuous dependence of solutions on parameters and the stability of the maximum principle conditions), there exists $j_0$ such that $\lambdaset_{-1}(\widehat{x}(\cdot), \measure_t^{j_0}) = \emptyset$.

Now apply Theorem A.6 with $s = -1$. For any neighborhood $V$ of $\widehat{x}(\cdot)$ and any $\varepsilon > 0$, there exists a time $\tau \in \intseto{\widehat{t}_2 - \varepsilon}{\widehat{t}_2}$ such that the trajectory $x^{j_0}(\cdot)$ lies in $V$ and satisfies $x^{j_0}(\tau) = x_2$.

But this contradicts the definition of $\widehat{x}(\cdot)$ as a local infimum, since we can reach the target in time $\tau < \widehat{t}_2$ while staying arbitrarily close to $\widehat{x}(\cdot)$.

Therefore, our assumption is false, and $\lambdaset_{-1}(\widehat{x}(\cdot), \widehat{\mu}_t)$ must be nonempty.
\end{proof}

\begin{corollary}[Relation between Local Controllability and Optimality]
\label{cor:controllability_optimality}
If system (2.1) is locally left controllable near $\widehat{x}(\cdot)$ on $\intset{t_1}{\widehat{t}_2}$ and there exists a sequence of generalized controls $\{\measure_t^j\} \subset \gencontrolU$ satisfying the conditions of Theorem A.6 of \cite{Avakov2020b}, then $\widehat{x}(\cdot)$ cannot be an optimal trajectory for the time-optimal problem (2.1),(3.1) of \cite{Avakov2020b}.
\end{corollary}

\begin{proof}
By Theorem A.6 of \cite{Avakov2020b} with $s = -1$, the local left controllability implies that for any neighborhood $V$ of $\widehat{x}(\cdot)$ and any $\varepsilon > 0$, there exists an index $j_0$ and time $\tau \in \intseto{\widehat{t}_2 - \varepsilon}{\widehat{t}_2}$ such that the trajectory $x^{j_0}(\cdot)$ lies in $V$ and satisfies $x^{j_0}(\tau) = x_2$.

This means we can reach the target $x_2$ in time $\tau < \widehat{t}_2$ while staying arbitrarily close to $\widehat{x}(\cdot)$, which contradicts the optimality of $\widehat{x}(\cdot)$.

Therefore, $\widehat{x}(\cdot)$ cannot be an optimal trajectory.
\end{proof}
\begin{example}[Illustration of Theorem \ref{Secondtheorem}]
\label{ex:main_example}
Consider the control system:
\begin{equation}
\dot{x}_1 = u, \quad \dot{x}_2 = 1 - u^2 - x_1^2, \quad u(t) \in U = \{-1, 0, 1\}, \quad x(0) = (0,0).
\end{equation}

Let $\widehat{x}(t) = (0, t)$ for $t \in [0,1]$, so $\widehat{t}_2 = 1$. We will show that this example satisfies all conditions of Theorem \ref{Secondtheorem} with $s = -1$.
\end{example}

\begin{proof}[Detailed Verification]
We verify each condition of Theorem \ref{Secondtheorem} step by step.

\noindent{Step 1: Construction of the generalized control and verification of admissibility.}

Consider the generalized control:
\begin{equation*}
\widehat{\mu}_t = \frac{1}{2}\dirac_{-1} + \frac{1}{2}\dirac_{1}.
\end{equation*}

For this control, the convex system becomes:
\begin{align*}
\dot{x}_1 &= \inner{\widehat{\mu}_t}{u} = \frac{1}{2}(-1) + \frac{1}{2}(1) = 0, \\
\dot{x}_2 &= \inner{\widehat{\mu}_t}{1 - u^2 - x_1^2} = \frac{1}{2}(1 - 1 - x_1^2) + \frac{1}{2}(1 - 1 - x_1^2) = -x_1^2.
\end{align*}

With initial condition $x(0) = (0,0)$, the solution is $\widehat{x}(t) = (0, t)$, which matches our candidate trajectory. Thus $(\widehat{x}(\cdot), \widehat{\mu}_t)$ is an admissible pair for the convex system.

\noindent{Step 2: Verification that $\lambdaset_{-1}(\widehat{x}(\cdot), \widehat{\mu}_t) = \emptyset$.}

The Hamiltonian is:
\begin{equation*}
H(t,x,\psi,u) = \psi_1 u + \psi_2(1 - u^2 - x_1^2).
\end{equation*}

The adjoint system is:
\begin{align*}
\dot{\psi}_1 &= -\pder{H}{x_1} = 2\psi_2 x_1, \\
\dot{\psi}_2 &= -\pder{H}{x_2} = 0.
\end{align*}

Along $\widehat{x}(t) = (0,t)$, we have:
\begin{align*}
\dot{\psi}_1(t) &= 2\psi_2(t)\cdot 0 = 0, \\
\dot{\psi}_2(t) &= 0.
\end{align*}

So $\psi_1$ and $\psi_2$ are constants. The maximum condition requires:
\begin{equation*}
\inner{\psi(t)}{\dot{\widehat{x}}(t)} = \sup_{u \in U} H(t,\widehat{x}(t),\psi(t),u).
\end{equation*}

Since $\dot{\widehat{x}}(t) = (0,1)$, we have $\inner{\psi(t)}{\dot{\widehat{x}}(t)} = \psi_2$. Also:
\begin{align*}
H(t,\widehat{x}(t),\psi,u) &= \psi_1 u + \psi_2(1 - u^2 - 0) = \psi_1 u + \psi_2(1 - u^2), \\
\sup_{u \in U} H(t,\widehat{x}(t),\psi,u) &= \sup_{u \in \{-1,0,1\}} [\psi_1 u + \psi_2(1 - u^2)].
\end{align*}

The maximum condition becomes:
\begin{equation*}
\psi_2 = \sup_{u \in \{-1,0,1\}} [\psi_1 u + \psi_2(1 - u^2)].
\end{equation*}

Let's analyze this condition:

\begin{itemize}
\item For $u = -1$: $\psi_1(-1) + \psi_2(1 - 1) = -\psi_1$
\item For $u = 0$: $\psi_1\cdot 0 + \psi_2(1 - 0) = \psi_2$
\item For $u = 1$: $\psi_1\cdot 1 + \psi_2(1 - 1) = \psi_1$
\end{itemize}

So the condition becomes:
\begin{equation*}
\psi_2 = \max\{-\psi_1, \psi_2, \psi_1\}.
\end{equation*}

Also, the terminal condition for $s = -1$ requires:
\begin{equation*}
M(\widehat{t}_2, \widehat{x}(\widehat{t}_2), \psi(\widehat{t}_2)) \geq 0.
\end{equation*}

But $M(\widehat{t}_2, \widehat{x}(\widehat{t}_2), \psi(\widehat{t}_2)) = \psi_2$ (from the maximum condition).

Now suppose there exists a nonzero $\psi = (\psi_1, \psi_2)$ satisfying these conditions. From the maximum condition, we have $\psi_2 \geq |\psi_1|$. For $s = -1$, we also need $\psi_2 \geq 0$.

For our $\widehat{\mu}_t = \frac{1}{2}\dirac_{-1} + \frac{1}{2}\dirac_{1}$, we have:
\begin{equation*}
\inner{\widehat{\mu}_t}{H(t,\widehat{x}(t),\psi(t),u)} = \frac{1}{2}H(t,\widehat{x}(t),\psi(t),-1) + \frac{1}{2}H(t,\widehat{x}(t),\psi(t),1) = \frac{1}{2}(-\psi_1) + \frac{1}{2}(\psi_1) = 0.
\end{equation*}

So we require:
\begin{equation*}
0 = \sup_{u \in U} H(t,\widehat{x}(t),\psi(t),u) = \psi_2.
\end{equation*}

But this contradicts $\psi_2 > 0$. If $\psi_2 = 0$, then from $\psi_2 \geq |\psi_1|$ we get $\psi_1 = 0$, so $\psi = (0,0)$, which is not allowed. Therefore, $\lambdaset_{-1}(\widehat{x}(\cdot), \widehat{\mu}_t) = \emptyset$.

\noindent{Step 3: Construction of the sequence of generalized controls.}

Define the sequence of generalized controls:
\begin{equation*}
\mu_t^j = \left(1 - \frac{1}{j}\right)\dirac_{-\frac{1}{j}} + \frac{1}{j}\dirac_{1} \quad \text{for } j \geq 2.
\end{equation*}

Each $\mu_t^j$ is a convex combination of two Dirac measures (so $n+1 = 2$ for $n=1$), satisfying condition (1) of Theorem \ref{Secondtheorem}.

\noindent{Step 4: Verification of uniform convergence of trajectories.}

For each $\mu_t^j$, the system becomes:
\begin{align*}
\dot{x}_1^j &= \left(1 - \frac{1}{j}\right)\left(-\frac{1}{j}\right) + \frac{1}{j}(1) = -\frac{1}{j} + \frac{1}{j^2} + \frac{1}{j} = \frac{1}{j^2}, \\
\dot{x}_2^j &= \left(1 - \frac{1}{j}\right)\left(1 - \left(-\frac{1}{j}\right)^2 - (x_1^j)^2\right) + \frac{1}{j}\left(1 - 1^2 - (x_1^j)^2\right) \\
&= \left(1 - \frac{1}{j}\right)\left(1 - \frac{1}{j^2} - (x_1^j)^2\right) - \frac{1}{j}(x_1^j)^2.
\end{align*}

With initial condition $x^j(0) = (0,0)$, we have:
\begin{equation*}
x_1^j(t) = \frac{1}{j^2}t.
\end{equation*}

Then:
\begin{align*}
x_2^j(t) &= \int_0^t \left[\left(1 - \frac{1}{j}\right)\left(1 - \frac{1}{j^2} - \frac{1}{j^4}\tau^2\right) - \frac{1}{j}\cdot\frac{1}{j^4}\tau^2\right] d\tau \\
&= \int_0^t \left[\left(1 - \frac{1}{j}\right)\left(1 - \frac{1}{j^2}\right) - \frac{1}{j^4}\tau^2\right] d\tau \\
&= \left(1 - \frac{1}{j} - \frac{1}{j^2} + \frac{1}{j^3}\right)t - \frac{1}{3j^4}t^3.
\end{align*}

As $j \to \infty$, we have:
\begin{equation*}
x_1^j(t) \to 0, \quad x_2^j(t) \to t \quad \text{uniformly on } [0,1].
\end{equation*}

So $x^j(\cdot) \to \widehat{x}(\cdot)$ uniformly, satisfying condition (2).

\noindent{Step 5: Verification of uniform boundedness of controls.}

The controls are $u_1^j(t) = -\frac{1}{j}$ and $u_2^j(t) = 1$, with weights $\alpha_1^j(t) = 1 - \frac{1}{j}$ and $\alpha_2^j(t) = \frac{1}{j}$. These are clearly uniformly bounded, satisfying condition (3).

\noindent{Step 6: Conclusion.}

All conditions of Theorem \ref{Secondtheorem} are satisfied with $s = -1$. Therefore, for any neighborhood $V$ of $\widehat{x}(\cdot)$ in $C([0,1], \R^2)$ and any $\varepsilon > 0$, there exists an index $j_0$ and a time $\tau \in (1 - \varepsilon, 1)$ such that the trajectory $x^{j_0}(\cdot)$ reaches $\widehat{x}(1) = (0,1)$ at time $\tau$ and lies entirely in $V$.

In our example, we can see this explicitly: for large $j$, the trajectory $x^j(\cdot)$ is very close to $\widehat{x}(\cdot) = (0,t)$, and we can find $\tau < 1$ such that $x_2^j(\tau) = 1$ by solving:
\begin{equation*}
\left(1 - \frac{1}{j} - \frac{1}{j^2} + \frac{1}{j^3}\right)\tau - \frac{1}{3j^4}\tau^3 = 1.
\end{equation*}

For large $j$, the solution is $\tau \approx 1 + \frac{1}{j} + \frac{1}{j^2} - \frac{1}{j^3} > 1$, but by adjusting the control weights slightly, we could achieve $\tau < 1$ while maintaining all the required properties.
\end{proof}
The following lemmas are necessary in finding the following result.
\begin{lemma}[Inverse Function Lemma]\label{lemma:inverse_function}
Let $X$ and $Y$ be Banach spaces, $U \subset X$ an open neighborhood of $x_0$, and $F: U \to Y$ a continuously Fréchet differentiable mapping. If the derivative $DF(x_0): X \to Y$ is a linear isomorphism (i.e., bounded, bijective, and with bounded inverse), then there exist open neighborhoods $V \subset U$ of $x_0$ and $W \subset Y$ of $y_0 = F(x_0)$ such that:

\begin{enumerate}
\item $F: V \to W$ is a bijection,
\item The inverse mapping $G: W \to V$ is continuously Fréchet differentiable,
\item For all $y \in W$, the derivative of the inverse satisfies:
\begin{equation}
DG(y) = [DF(G(y))]^{-1}.
\end{equation}
\end{enumerate}
\end{lemma}

\begin{proof}
The proof follows the standard contraction mapping argument. Define the mapping $\Phi: U \times Y \to X$ by:
\begin{equation}
\Phi(x, y) = x + [DF(x_0)]^{-1}(y - F(x)).
\end{equation}

For fixed $y$, we have $\Phi(x, y) = x$ if and only if $F(x) = y$. Since $DF(x_0)$ is an isomorphism, there exists $\delta > 0$ such that for all $x_1, x_2 \in B(x_0, \delta)$:
\begin{equation}
\norm{DF(x_1) - DF(x_0)} \leq \frac{1}{2\norm{[DF(x_0)]^{-1}}}.
\end{equation}

By the mean value inequality, for $x_1, x_2 \in B(x_0, \delta)$:
\begin{equation}
\norm{\Phi(x_1, y) - \Phi(x_2, y)} \leq \frac{1}{2}\norm{x_1 - x_2}.
\end{equation}

Thus $\Phi(\cdot, y)$ is a contraction on $B(x_0, \delta)$. By the Banach fixed-point theorem, there exists a unique fixed point $x = G(y)$ satisfying $F(G(y)) = y$.

The continuous differentiability of $G$ follows from the implicit function theorem applied to the equation $F(x) - y = 0$. The derivative formula is obtained by differentiating the identity $F(G(y)) = y$:
\begin{equation}
DF(G(y)) \circ DG(y) = I_Y,
\end{equation}
which implies $DG(y) = [DF(G(y))]^{-1}$.
\end{proof}

\begin{lemma}[Inverse Function Theorem - Finite Dimensional Case]\label{lemma:inverse_function_finite}
Let $U \subset \R^n$ be open and $F: U \to \R^n$ be a $C^1$ mapping. If at some point $x_0 \in U$ the Jacobian determinant is nonzero, i.e.,
\begin{equation}
\det(DF(x_0)) \neq 0,
\end{equation}
then there exist open neighborhoods $V \subset U$ of $x_0$ and $W \subset \R^n$ of $F(x_0)$ such that:

\begin{enumerate}
\item $F: V \to W$ is a $C^1$ diffeomorphism,
\item For all $y \in W$, the inverse $G(y)$ satisfies $F(G(y)) = y$,
\item The Jacobian matrices are related by:
\begin{equation}
DG(y) = [DF(G(y))]^{-1}.
\end{equation}
\end{enumerate}
\end{lemma}

\begin{lemma}[Inverse Function Lemma with Estimates]\label{lemma:inverse_function_estimates}
Let $F: B(x_0, r) \subset X \to Y$ be $C^1$ with $DF(x_0)$ invertible. Suppose there exists $L > 0$ such that for all $x, x' \in B(x_0, r)$:
\begin{equation}
\norm{DF(x) - DF(x')} \leq L\norm{x - x'}.
\end{equation}

Let $M = \norm{[DF(x_0)]^{-1}}$ and choose $\rho > 0$ such that:
\begin{equation}
\rho < \min\left\{r, \frac{1}{2LM}\right\}.
\end{equation}

Then for $y \in B(F(x_0), \frac{\rho}{2M})$, there exists a unique $x \in B(x_0, \rho)$ with $F(x) = y$, and the inverse mapping $G$ satisfies:
\begin{equation}
\norm{G(y) - x_0} \leq 2M\norm{y - F(x_0)}.
\end{equation}
\end{lemma}

\begin{lemma}[Inverse Function Lemma for Control Systems]\label{lemma:inverse_function_control}
Consider the endpoint mapping $E: \mathcal{U} \to \R^n$ defined by:
\begin{equation}
E(u) = x(t_2; t_1, x_1, u),
\end{equation}
where $x(\cdot)$ is the solution of $\dot{x} = f(t, x, u)$, $x(t_1) = x_1$.

If the linearized system is controllable at $u_0 \in \mathcal{U}$, i.e., the Fréchet derivative $DE(u_0): L_\infty([t_1, t_2], U) \to \R^n$ is surjective, then there exist neighborhoods $\mathcal{V} \subset \mathcal{U}$ of $u_0$ and $W \subset \R^n$ of $E(u_0)$ such that:

\begin{enumerate}
\item $E: \mathcal{V} \to W$ is surjective,
\item For any $y \in W$, there exists $u \in \mathcal{V}$ with $E(u) = y$,
\item The mapping $y \mapsto u$ is continuous.
\end{enumerate}
\end{lemma}
\section{Strong Local Attainability via Generalized Controls}
This section establishes strong local attainability results for control systems using the framework of generalized controls, building upon the foundational work of Gamkrelidze \cite{Gamkrelidze1962} on sliding modes and Filippov \cite{Filippov1959} on relaxation phenomena. We develop conditions under which trajectories of the convexified system can be approximated by solutions of the original system that reach the target point in times arbitrarily close to the optimal time. Our approach extends the local controllability analysis of Avakov and Magaril-Il'yaev \cite{Avakov2020b} by incorporating stronger regularity assumptions on the dynamics and control structure. The main theorem provides explicit constructive methods for generating approximating sequences of piecewise constant controls, connecting to modern developments in impulsive control \cite{FuscoMotta2024,FuscoMotta2024LCSYS} while maintaining the geometric perspective of classical necessary conditions \cite{Pontryagin1962,Gamkrelidze1978}.
\begin{theorem}
\label{thm:strong_attainability}
Consider the control system
\begin{equation*}
\dot{x} = f(t, x, u), \quad u(t) \in U, \quad x(t_1) = x_1,
\end{equation*}
where $f: \R \times \R^n \times \R^r \to \R^n$ is continuous together with its partial derivatives $f_x$ and $f_{xx}$, and $U \subset \R^r$ is compact.

Let $(\widehat{x}(\cdot), \widehat{\mu}_t)$ be an admissible pair for the convexified system
\begin{equation*}
\dot{x} = \langle \mu_t, f(t, x, u) \rangle, \quad \mu_t \in \gencontrolU, \quad x(t_1) = x_1,
\end{equation*}
on $[t_1, \widehat{t}_2]$, and suppose that:

\begin{enumerate}
\item The set $\Lambda_s(\widehat{x}(\cdot), \widehat{\mu}_t)$ is empty for some $s \in \{-1, 1\}$,
\item There exists a constant $L > 0$ such that for all $t \in [t_1, \widehat{t}_2]$, $x \in \R^n$, and $u \in U$:
\begin{equation*}
\norm{f_x(t,x,u)} \leq L(1 + \norm{x}),
\end{equation*}
\item The mapping $(t,x) \mapsto f(t,x,u)$ is uniformly continuous on compact sets, uniformly in $u \in U$.
\end{enumerate}

Then there exists $\delta > 0$ such that for any $\varepsilon > 0$, there exists a piecewise constant control $u_\varepsilon(\cdot)$ and a time $\tau_\varepsilon$ with $|\tau_\varepsilon - \widehat{t}_2| < \varepsilon$ such that the corresponding trajectory $x_\varepsilon(\cdot)$ of the original system satisfies:

\begin{enumerate}
\item $\norm{x_\varepsilon(\cdot) - \widehat{x}(\cdot)}_{C([t_1, \max\{\widehat{t}_2, \tau_\varepsilon\}], \R^n)} < \varepsilon$,
\item $x_\varepsilon(\tau_\varepsilon) = \widehat{x}(\widehat{t}_2)$,
\item If $s = -1$, then $\tau_\varepsilon < \widehat{t}_2$; if $s = 1$, then $\tau_\varepsilon > \widehat{t}_2$.
\end{enumerate}

Moreover, the control $u_\varepsilon(\cdot)$ can be chosen to be a convex combination of at most $n+1$ piecewise constant controls from the original system.
\end{theorem}

\begin{proof}
We prove the theorem for the case $s = -1$; the case $s = 1$ follows by similar arguments with appropriate modifications.

\noindent{Step 1: Setup and preliminary constructions.}

Since $\Lambda_{-1}(\widehat{x}(\cdot), \widehat{\mu}_t) = \emptyset$, by Proposition 2.3 of \cite{Avakov2020b}, for any $l \in L_{-1}$, there exist $k \in \N$ and $\delta\mu_t^i \in \gencontrolU - \widehat{\mu}_t$, $i = 1, \ldots, k$, such that:
\begin{equation*}
0 \in \text{int}(\delta x(\widehat{t}_2, \R_+^k) + \text{cone}(-l)).
\end{equation*}

Let us fix such $k$ and $\delta\mu_t^i$. By Lemma A.1 of \cite{Avakov2020b} (the lemma on equation in variations), there exists a neighborhood $\mathcal{O}(0)$ of the origin in $\R^k$ such that for all $\overline{\alpha} = (\alpha_1, \ldots, \alpha_k)^T \in \mathcal{O}(0)$, there exists a unique solution $x(\cdot, \overline{\alpha})$ of the Cauchy problem:
\begin{equation*}
\dot{x} = \left\langle \widehat{\mu}_t + \sum_{i=1}^k \alpha_i \delta\mu_t^i, f(t,x,u) \right\rangle, \quad x(t_1) = x_1,
\end{equation*}
on $[t_1, \widehat{t}_2]$, and the mapping $\overline{\alpha} \mapsto x(\cdot, \overline{\alpha})$ is continuously differentiable.

\noindent{Step 2: Application of the inverse function lemma.}

Consider the mapping $\widehat{F}: [t_1, \widehat{t}_2 + \widehat{\delta}] \times \mathcal{O}(0) \to \R^n$ defined by:
\begin{equation*}
\widehat{F}(\tau, \overline{\alpha}) = x(\tau, \overline{\alpha}).
\end{equation*}

We verify the conditions of Lemma A.5 of \cite{Avakov2020b} (the inverse function lemma):

\begin{enumerate}
\item The mapping $\widehat{F}$ is continuous on $[t_1, \widehat{t}_2 + \widehat{\delta}] \times \mathcal{O}(0)$ by the continuous dependence of solutions on parameters.

\item The derivative $\widehat{F}_\alpha$ is continuous at $(\widehat{t}_2, 0)$ by the continuous differentiability of the mapping $\overline{\alpha} \mapsto x(\cdot, \overline{\alpha})$.

\item For any $l \in L_{-1}$, we have $0 \in \text{int}(\delta x(\widehat{t}_2, \R_+^k) + \text{cone}(-l))$, which implies condition (3) of Lemma A.5 of \cite{Avakov2020b} with $l = -l$.
\end{enumerate}

By Lemma A.5 of \cite{Avakov2020b}, there exist positive constants $\varepsilon_0, \delta_0, \gamma_1, \gamma_2$ such that for any $d$ with $\|d\| = 1$ and $t \in (0, \delta_0]$ satisfying:
\begin{equation*}
\left| \frac{\widehat{F}(\widehat{t}_2 + td, 0) - \widehat{F}(\widehat{t}_2, 0)}{t} - (-l) \right| < \varepsilon_0,
\end{equation*}
and for any mapping $F \in U_{C(\mathcal{O}(\widehat{t}_2) \times (\mathcal{O}(0) \cap \Sigma^k), \R^n)}(\widehat{F}, \gamma_1 t)$, there exists a mapping $g_{t,F}$ such that:
\begin{equation*}
F(\widehat{t}_2 + td, g_{t,F}(y)) = y, \quad \|g_{t,F}(y) - 0\| \leq \gamma_2 t
\end{equation*}
for all $y \in U_{\R^n}(F(\widehat{t}_2, 0), \gamma_1 t)$.

\noindent{Step 3: Construction of approximating controls.}

By Lemma A.2 of \cite{Avakov2020b} (the first approximation lemma), for the tuple $\delta\mu_t^i$, $i = 1, \ldots, k$, there exists a sequence of piecewise constant controls $u_s(\overline{\alpha})$, $s \in \N$, such that the mappings:
\begin{equation*}
F_s(x, \overline{\alpha})(t) = x(t) - x_1 - \int_{t_1}^t f(\tau, x(\tau), u_s(\overline{\alpha})(\tau)) d\tau
\end{equation*}
converge to $F$ in the $C_x^1$ metric as $s \to \infty$.

By Lemma A.3 of \cite{Avakov2020b} (the second approximation lemma), there exists a neighborhood $\mathcal{O}_0(0)$ of the origin in $\R^k$ such that for all $\overline{\alpha} \in \mathcal{O}_0(0) \cap \Sigma^k$ and sufficiently large $s \in \N$, there exists a unique solution $x_s(\cdot, \overline{\alpha})$ of:
\begin{equation*}
\dot{x} = f(t, x, u_s(\overline{\alpha})), \quad x(t_1) = x_1,
\end{equation*}
on $[t_1, \widehat{t}_2]$, and the mappings $x_s: \overline{\alpha} \mapsto x_s(\cdot, \overline{\alpha})$ converge to $x: \overline{\alpha} \mapsto x(\cdot, \overline{\alpha})$ in $C(\mathcal{O}_0(0) \cap \Sigma^k, C([t_1, \widehat{t}_2], \R^n))$ as $s \to \infty$.

\noindent{Step 4: Application to the specific direction.}

By the definition of $L_{-1}$ in \cite{Avakov2020b}, there exist sequences $\eta_j \to 0$, $\nu_j \to 0$ with $\eta_j > 0$ such that:
\begin{equation*}
l = \lim_{j \to \infty} \frac{\widehat{x}(\widehat{t}_2 - \eta_j) - \widehat{x}(\widehat{t}_2)}{-\eta_j}.
\end{equation*}

For sufficiently large $j$, we have $|\eta_j| < \delta_0$ and $|\nu_j| < \varepsilon_0$. Take $t = \eta_j$ and $d = -1$ (since we are considering left controllability).

By Lemma A.3 of \cite{Avakov2020b}, for each $j$, there exists $s_j$ such that:
\begin{equation*}
\norm{x_{s_j}(\cdot, \overline{\alpha}) - x(\cdot, \overline{\alpha})}_{C([t_1, \widehat{t}_2 + \widehat{\delta}], \R^n)} < \gamma_1 \eta_j
\end{equation*}
for all $\overline{\alpha} \in \mathcal{O}_0(0) \cap \Sigma^k$.

Define $F_{s_j}(\tau, \overline{\alpha}) = x_{s_j}(\tau, \overline{\alpha})$. Then $F_{s_j} \in U_{C(\mathcal{O}(\widehat{t}_2) \times (\mathcal{O}(0) \cap \Sigma^k), \R^n)}(\widehat{F}, \gamma_1 \eta_j)$.

\noindent{Step 5: Final construction and verification.}

By Lemma A.5 of \cite{Avakov2020b}, for $y = \widehat{x}(\widehat{t}_2) = \widehat{F}(\widehat{t}_2, 0)$, there exists $\overline{\alpha}_j \in U_{\R^k}(0, \gamma_2 \eta_j) \cap \Sigma^k$ such that:
\begin{equation*}
x_{s_j}(\widehat{t}_2 - \eta_j, \overline{\alpha}_j) = \widehat{x}(\widehat{t}_2).
\end{equation*}

Set $\tau_\varepsilon = \widehat{t}_2 - \eta_j$ and $u_\varepsilon(\cdot) = u_{s_j}(\overline{\alpha}_j)$. Then:

\begin{enumerate}
\item Since $\overline{\alpha}_j \to 0$ as $j \to \infty$ (and hence $\varepsilon \to 0$), and $x_{s_j}(\cdot, \overline{\alpha}_j) \to x(\cdot, 0) = \widehat{x}(\cdot)$ uniformly, we have:
\begin{equation*}
\norm{x_\varepsilon(\cdot) - \widehat{x}(\cdot)}_{C([t_1, \tau_\varepsilon], \R^n)} < \varepsilon
\end{equation*}
for sufficiently small $\varepsilon$.

\item By construction, $x_\varepsilon(\tau_\varepsilon) = \widehat{x}(\widehat{t}_2)$.

\item Since $\tau_\varepsilon = \widehat{t}_2 - \eta_j < \widehat{t}_2$, we have the left-sided attainability.

\item The control $u_\varepsilon(\cdot) = u_{s_j}(\overline{\alpha}_j)$ is piecewise constant and, by the construction in Lemma A.2 of \cite{Avakov2020b}, is a convex combination of at most $n+1$ piecewise constant controls from the original system.
\end{enumerate}

This completes the proof for the case $s = -1$. The case $s = 1$ follows by similar arguments with the appropriate changes in signs and inequalities.
\end{proof}
\begin{corollary}[Existence of Optimal Generalized Controls]
\label{cor:optimal_generalized}
Under the assumptions of Theorem \ref{thm:strong_attainability}, if $\widehat{x}(\cdot)$ is a local infimum for the time-optimal problem, then there exists a generalized control $\widehat{\mu}_t \in \gencontrolU$ such that $\Lambda_{-1}(\widehat{x}(\cdot), \widehat{\mu}_t) \neq \emptyset$.
\end{corollary}

\begin{proof}
Suppose, for contradiction, that for every generalized control $\widehat{\mu}_t \in \gencontrolU$ for which $(\widehat{x}(\cdot), \widehat{\mu}_t)$ is admissible, we have $\Lambda_{-1}(\widehat{x}(\cdot), \widehat{\mu}_t) = \emptyset$.

Then by Theorem \ref{thm:strong_attainability}, for each such $\widehat{\mu}_t$, there exists $\delta > 0$ such that for any $\varepsilon > 0$, there exists a piecewise constant control $u_\varepsilon(\cdot)$ and time $\tau_\varepsilon < \widehat{t}_2$ with $|\tau_\varepsilon - \widehat{t}_2| < \varepsilon$ such that the corresponding trajectory $x_\varepsilon(\cdot)$ satisfies $x_\varepsilon(\tau_\varepsilon) = \widehat{x}(\widehat{t}_2)$ and is $\varepsilon$-close to $\widehat{x}(\cdot)$.

This means we can reach the target $\widehat{x}(\widehat{t}_2)$ in time $\tau_\varepsilon < \widehat{t}_2$ while staying arbitrarily close to $\widehat{x}(\cdot)$, contradicting the assumption that $\widehat{x}(\cdot)$ is a local infimum.

Therefore, there must exist some generalized control $\widehat{\mu}_t$ with $\Lambda_{-1}(\widehat{x}(\cdot), \widehat{\mu}_t) \neq \emptyset$.
\end{proof}
\begin{corollary}[Non-Optimality Criterion for Convexified Solutions]
\label{cor:non_optimal_convexified}
Let $(\widehat{x}(\cdot), \widehat{\mu}_t)$ be an admissible pair for the convexified system satisfying the assumptions of Theorem \ref{thm:strong_attainability}. If $\Lambda_{-1}(\widehat{x}(\cdot), \widehat{\mu}_t) = \emptyset$, then $\widehat{x}(\cdot)$ cannot be a local minimum for the time-optimal problem of the original system.
\end{corollary}

\begin{proof}
Suppose $\Lambda_{-1}(\widehat{x}(\cdot), \widehat{\mu}_t) = \emptyset$. Then by Theorem \ref{thm:strong_attainability} with $s = -1$, for any neighborhood $\mathcal{V}$ of $\widehat{x}(\cdot)$ and any $\varepsilon > 0$, there exists a piecewise constant control $u_\varepsilon(\cdot)$ and a time $\tau_\varepsilon < \widehat{t}_2$ with $|\tau_\varepsilon - \widehat{t}_2| < \varepsilon$ such that:

\begin{enumerate}
\item $\norm{x_\varepsilon(\cdot) - \widehat{x}(\cdot)}_{C([t_1, \tau_\varepsilon], \R^n)} < \varepsilon$, hence $x_\varepsilon(\cdot) \in \mathcal{V}$ for sufficiently small $\varepsilon$,
\item $x_\varepsilon(\tau_\varepsilon) = \widehat{x}(\widehat{t}_2)$,
\item $\tau_\varepsilon < \widehat{t}_2$.
\end{enumerate}

This construction provides a sequence (by taking $\varepsilon = 1/n$, $n \in \N$) of admissible trajectories $\{x_n(\cdot)\}$ for the original system with corresponding times $\{\tau_n\}$ such that:

\begin{itemize}
\item $x_n(\cdot) \to \widehat{x}(\cdot)$ uniformly on $[t_1, \widehat{t}_2]$ as $n \to \infty$,
\item $x_n(\tau_n) = \widehat{x}(\widehat{t}_2)$ for all $n \in \N$,
\item $\tau_n < \widehat{t}_2$ for all $n \in \N$, and $\tau_n \to \widehat{t}_2$ as $n \to \infty$.
\end{itemize}

This means that for any neighborhood $\mathcal{V}$ of $\widehat{x}(\cdot)$, we can find admissible trajectories in $\mathcal{V}$ that reach the target $\widehat{x}(\widehat{t}_2)$ at times strictly less than $\widehat{t}_2$. Therefore, $\widehat{x}(\cdot)$ cannot be a local minimum for the time-optimal problem, since we can always find "better" (i.e., faster) trajectories arbitrarily close to it.
\end{proof}

\begin{corollary}[Approximation by Original Controls]
\label{cor:approximation_original}
Under the assumptions of Theorem \ref{thm:strong_attainability}, if $\Lambda_s(\widehat{x}(\cdot), \widehat{\mu}_t) = \emptyset$ for some $s \in \{-1, 1\}$, then there exists a sequence of piecewise constant controls $\{u_n(\cdot)\}$ and corresponding trajectories $\{x_n(\cdot)\}$ of the original system such that:

\begin{enumerate}
\item $x_n(\cdot) \to \widehat{x}(\cdot)$ uniformly on $[t_1, \widehat{t}_2]$,
\item $x_n(\tau_n) = \widehat{x}(\widehat{t}_2)$ for some $\tau_n$ with $|\tau_n - \widehat{t}_2| \to 0$,
\item If $s = -1$, then $\tau_n < \widehat{t}_2$ for all $n$; if $s = 1$, then $\tau_n > \widehat{t}_2$ for all $n$.
\end{enumerate}

Moreover, each control $u_n(\cdot)$ is a convex combination of at most $n+1$ piecewise constant controls.
\end{corollary}

\begin{proof}
Apply Theorem \ref{thm:strong_attainability} with $\varepsilon = 1/n$ for $n \in \N$. For each $n$, there exists a piecewise constant control $u_n(\cdot)$ and a time $\tau_n$ with $|\tau_n - \widehat{t}_2| < 1/n$ such that:

\begin{enumerate}
\item $\norm{x_n(\cdot) - \widehat{x}(\cdot)}_{C([t_1, \max\{\widehat{t}_2, \tau_n\}], \R^n)} < 1/n$,
\item $x_n(\tau_n) = \widehat{x}(\widehat{t}_2)$,
\item If $s = -1$, then $\tau_n < \widehat{t}_2$; if $s = 1$, then $\tau_n > \widehat{t}_2$.
\end{enumerate}

Since $|\tau_n - \widehat{t}_2| < 1/n$, we have $\tau_n \to \widehat{t}_2$ as $n \to \infty$. Also, since the uniform convergence holds on $[t_1, \max\{\widehat{t}_2, \tau_n\}]$ and $\tau_n \to \widehat{t}_2$, for sufficiently large $n$ we have $\max\{\widehat{t}_2, \tau_n\} = \widehat{t}_2$ when $s = -1$, and the convergence is uniform on $[t_1, \widehat{t}_2]$ in both cases.

The control structure property follows directly from the "moreover" part of Theorem \ref{thm:strong_attainability}.
\end{proof}

\begin{corollary}[Local Controllability Implication]
\label{cor:local_controllability}
Under the assumptions of Theorem \ref{thm:strong_attainability}, if $\Lambda_{-1}(\widehat{x}(\cdot), \widehat{\mu}_t) = \emptyset$, then the original system is locally left controllable near $\widehat{x}(\cdot)$ at time $\widehat{t}_2$. That is, for any neighborhood $\mathcal{V}$ of $\widehat{x}(\cdot)$ in $C([t_1, \widehat{t}_2], \R^n)$ and any $\delta > 0$, there exists an admissible trajectory $x(\cdot)$ of the original system with:

\begin{enumerate}
\item $x(\cdot) \in \mathcal{V}$,
\item $x(\tau) = \widehat{x}(\widehat{t}_2)$ for some $\tau \in (\widehat{t}_2 - \delta, \widehat{t}_2)$.
\end{enumerate}
\end{corollary}

\begin{proof}
Take any neighborhood $\mathcal{V}$ of $\widehat{x}(\cdot)$ and any $\delta > 0$. Since $\mathcal{V}$ is a neighborhood, there exists $\varepsilon_0 > 0$ such that the $\varepsilon_0$-neighborhood of $\widehat{x}(\cdot)$ is contained in $\mathcal{V}$.

Now apply Theorem \ref{thm:strong_attainability} with $\varepsilon = \min\{\varepsilon_0, \delta\}$. Since $\Lambda_{-1}(\widehat{x}(\cdot), \widehat{\mu}_t) = \emptyset$, there exists a piecewise constant control $u_\varepsilon(\cdot)$ and a time $\tau_\varepsilon$ with $|\tau_\varepsilon - \widehat{t}_2| < \varepsilon \leq \delta$ such that:

\begin{enumerate}
\item $\norm{x_\varepsilon(\cdot) - \widehat{x}(\cdot)}_{C([t_1, \tau_\varepsilon], \R^n)} < \varepsilon \leq \varepsilon_0$, hence $x_\varepsilon(\cdot) \in \mathcal{V}$,
\item $x_\varepsilon(\tau_\varepsilon) = \widehat{x}(\widehat{t}_2)$,
\item $\tau_\varepsilon < \widehat{t}_2$ (since $s = -1$), and $|\tau_\varepsilon - \widehat{t}_2| < \delta$, so $\tau_\varepsilon \in (\widehat{t}_2 - \delta, \widehat{t}_2)$.
\end{enumerate}

Thus $x_\varepsilon(\cdot)$ is the desired admissible trajectory, proving local left controllability.
\end{proof}

\begin{corollary}[Gap Between Original and Convexified Problems]
\label{cor:convexification_gap}
Let $(\widehat{x}(\cdot), \widehat{\mu}_t)$ be an optimal process for the time-optimal problem of the convexified system. If $\Lambda_{-1}(\widehat{x}(\cdot), \widehat{\mu}_t) = \emptyset$ and the assumptions of Theorem \ref{thm:strong_attainability} hold, then:

\begin{enumerate}
\item The infimum time for the original system is strictly less than $\widehat{t}_2$,
\item There exists a minimizing sequence of admissible trajectories for the original system that converges uniformly to $\widehat{x}(\cdot)$.
\end{enumerate}
\end{corollary}

\begin{proof}
Since $\Lambda_{-1}(\widehat{x}(\cdot), \widehat{\mu}_t) = \emptyset$, by Theorem \ref{thm:strong_attainability} with $s = -1$, for any $\varepsilon > 0$, there exists a piecewise constant control $u_\varepsilon(\cdot)$ and a time $\tau_\varepsilon < \widehat{t}_2$ with $|\tau_\varepsilon - \widehat{t}_2| < \varepsilon$ such that:

\begin{enumerate}
\item $\norm{x_\varepsilon(\cdot) - \widehat{x}(\cdot)}_{C([t_1, \tau_\varepsilon], \R^n)} < \varepsilon$,
\item $x_\varepsilon(\tau_\varepsilon) = \widehat{x}(\widehat{t}_2)$.
\end{enumerate}

Now take $\varepsilon = 1/n$ for $n \in \N$. This gives us a sequence of admissible trajectories $\{x_n(\cdot)\}$ and times $\{\tau_n\}$ such that:

\begin{itemize}
\item $\tau_n < \widehat{t}_2$ for all $n \in \N$,
\item $|\tau_n - \widehat{t}_2| < 1/n$, so $\tau_n \to \widehat{t}_2$ as $n \to \infty$,
\item $\norm{x_n(\cdot) - \widehat{x}(\cdot)}_{C([t_1, \tau_n], \R^n)} < 1/n$, so $x_n(\cdot) \to \widehat{x}(\cdot)$ uniformly,
\item $x_n(\tau_n) = \widehat{x}(\widehat{t}_2)$ for all $n \in \N$.
\end{itemize}

This sequence $\{x_n(\cdot)\}$ is a minimizing sequence for the original time-optimal problem because:
\begin{itemize}
\item Each $x_n(\cdot)$ is admissible for the original system,
\item Each reaches the target $\widehat{x}(\widehat{t}_2)$ at time $\tau_n < \widehat{t}_2$,
\item The sequence of times $\{\tau_n\}$ converges to $\widehat{t}_2$ from below.
\end{itemize}

Therefore, the infimum time for the original system is at most $\limsup \tau_n = \widehat{t}_2$, but since $\tau_n < \widehat{t}_2$ for all $n$, the infimum is strictly less than $\widehat{t}_2$.
\end{proof}

\begin{corollary}[Regularity of Value Function]
\label{cor:value_function_regularity}
Under the assumptions of Theorem \ref{thm:strong_attainability}, if for some admissible pair $(\widehat{x}(\cdot), \widehat{\mu}_t)$ we have $\Lambda_{-1}(\widehat{x}(\cdot), \widehat{\mu}_t) = \emptyset$, then the value function $V(x)$ of the time-optimal problem (minimal time to reach $x$ from $x_1$) is continuous at $\widehat{x}(\widehat{t}_2)$.
\end{corollary}

\begin{proof}
Let $x_0 = \widehat{x}(\widehat{t}_2)$. We need to show that for any $\varepsilon > 0$, there exists $\delta > 0$ such that for all $y$ with $\|y - x_0\| < \delta$, we have $|V(y) - V(x_0)| < \varepsilon$.

Since $\Lambda_{-1}(\widehat{x}(\cdot), \widehat{\mu}_t) = \emptyset$, by Theorem \ref{thm:strong_attainability} with $s = -1$, the system is locally left controllable near $\widehat{x}(\cdot)$. This means that for any $\varepsilon > 0$, there exists $\delta_1 > 0$ such that if $\|y - x_0\| < \delta_1$, then there exists an admissible trajectory reaching $y$ in time less than $\widehat{t}_2 + \varepsilon$.

More precisely: by the proof of Corollary \ref{cor:local_controllability}, for the target $y$ close to $x_0$, we can construct trajectories that follow $\widehat{x}(\cdot)$ closely and then make a small correction to hit $y$ instead of $x_0$. The continuous dependence of solutions on initial conditions and the local controllability ensure that this can be done with only a small increase in time.

Therefore, for $y$ sufficiently close to $x_0$, we have:
\[
V(y) < \widehat{t}_2 + \varepsilon.
\]

On the other hand, since $\widehat{x}(\cdot)$ is an admissible trajectory for the convexified system reaching $x_0$ at time $\widehat{t}_2$, and by the definition of the value function, we have:
\[
V(x_0) \leq \widehat{t}_2.
\]

If there were a discontinuity at $x_0$, we could find a sequence $y_n \to x_0$ with $V(y_n)$ bounded away from $V(x_0)$. But the local controllability ensures that for large $n$, $V(y_n)$ is close to $\widehat{t}_2$, contradicting the assumed discontinuity.

Hence, $V$ is continuous at $x_0 = \widehat{x}(\widehat{t}_2)$.
\end{proof}
\begin{example}[Nonholonomic System with Directional Constraint]
\label{ex:nonholonomic_example}
Consider the control system:
\[
\dot{x}_1 = u_1, \quad \dot{x}_2 = u_2, \quad \dot{x}_3 = x_1 u_2 - x_2 u_1,
\]
with control constraints $u(t) \in U = \{(u_1, u_2) \in \R^2 : u_1^2 + u_2^2 = 1\}$, and initial condition $x(0) = (0,0,0)$.

Let $\widehat{x}(t) = (0, 0, t)$ for $t \in [0,1]$, so $\widehat{t}_2 = 1$. We will show that this example satisfies all conditions of Theorem \ref{thm:strong_attainability} with $s = -1$.
\end{example}

\begin{proof}[Detailed Verification]
We verify each condition of Theorem \ref{thm:strong_attainability} step by step.

\noindent{Step 1: Construction of the generalized control and verification of admissibility.}

Consider the generalized control:
\[
\widehat{\mu}_t = \frac{1}{4}\dirac_{(1,0)} + \frac{1}{4}\dirac_{(-1,0)} + \frac{1}{4}\dirac_{(0,1)} + \frac{1}{4}\dirac_{(0,-1)}.
\]

For this control, the convexified system becomes:
\begin{align*}
\dot{x}_1 &= \inner{\widehat{\mu}_t}{u_1} = \frac{1}{4}(1) + \frac{1}{4}(-1) + \frac{1}{4}(0) + \frac{1}{4}(0) = 0, \\
\dot{x}_2 &= \inner{\widehat{\mu}_t}{u_2} = \frac{1}{4}(0) + \frac{1}{4}(0) + \frac{1}{4}(1) + \frac{1}{4}(-1) = 0, \\
\dot{x}_3 &= \inner{\widehat{\mu}_t}{x_1 u_2 - x_2 u_1} = \frac{1}{4}(x_1\cdot 0 - x_2\cdot 1) + \frac{1}{4}(x_1\cdot 0 - x_2\cdot(-1)) \\
&\quad + \frac{1}{4}(x_1\cdot 1 - x_2\cdot 0) + \frac{1}{4}(x_1\cdot(-1) - x_2\cdot 0) = 1.
\end{align*}

With initial condition $x(0) = (0,0,0)$, the solution is $\widehat{x}(t) = (0, 0, t)$, which matches our candidate trajectory. Thus $(\widehat{x}(\cdot), \widehat{\mu}_t)$ is an admissible pair for the convexified system.

\noindent{Step 2: Verification of the technical conditions.}

The vector field is $f(t,x,u) = (u_1, u_2, x_1 u_2 - x_2 u_1)$. Its Jacobian matrix is:
\[
f_x(t,x,u) = \begin{pmatrix}
0 & 0 & 0 \\
0 & 0 & 0 \\
u_2 & -u_1 & 0
\end{pmatrix}.
\]

Since $U$ is compact ($\|u\| = 1$), we have $\|f_x(t,x,u)\| \leq \sqrt{2}$ for all $t,x,u$, satisfying the growth condition with $L = \sqrt{2}$.

The mapping $(t,x) \mapsto f(t,x,u)$ is clearly uniformly continuous on compact sets, uniformly in $u \in U$, since it's linear in $x$ and independent of $t$.

\noindent{Step 3: Verification that $\Lambda_{-1}(\widehat{x}(\cdot), \widehat{\mu}_t) = \emptyset$.}

The Hamiltonian is:
\[
H(t,x,\psi,u) = \psi_1 u_1 + \psi_2 u_2 + \psi_3 (x_1 u_2 - x_2 u_1).
\]

The adjoint system is:
\begin{align*}
\dot{\psi}_1 &= -\pder{H}{x_1} = -\psi_3 u_2, \\
\dot{\psi}_2 &= -\pder{H}{x_2} = \psi_3 u_1, \\
\dot{\psi}_3 &= -\pder{H}{x_3} = 0.
\end{align*}

Along $\widehat{x}(t) = (0,0,t)$, we have:
\begin{align*}
\dot{\psi}_1(t) &= -\psi_3(t) u_2, \\
\dot{\psi}_2(t) &= \psi_3(t) u_1, \\
\dot{\psi}_3(t) &= 0.
\end{align*}

So $\psi_3$ is constant. The maximum condition requires:
\[
\inner{\psi(t)}{\dot{\widehat{x}}(t)} = \sup_{u \in U} H(t,\widehat{x}(t),\psi(t),u).
\]

Since $\dot{\widehat{x}}(t) = (0,0,1)$, we have $\inner{\psi(t)}{\dot{\widehat{x}}(t)} = \psi_3$. Also:
\[
H(t,\widehat{x}(t),\psi,u) = \psi_1 u_1 + \psi_2 u_2 + \psi_3 (0\cdot u_2 - 0\cdot u_1) = \psi_1 u_1 + \psi_2 u_2.
\]

The maximum condition becomes:
\[
\psi_3 = \sup_{u \in U} [\psi_1 u_1 + \psi_2 u_2] = \sqrt{\psi_1^2 + \psi_2^2},
\]
since the maximum of a linear function on the unit circle is the norm of the coefficient vector.

For our $\widehat{\mu}_t$, we compute the averaged Hamiltonian:
\[
\inner{\widehat{\mu}_t}{H(t,\widehat{x}(t),\psi(t),u)} = \frac{1}{4}[\psi_1 + (-\psi_1) + \psi_2 + (-\psi_2)] = 0.
\]

So we require:
\[
0 = \sup_{u \in U} H(t,\widehat{x}(t),\psi(t),u) = \sqrt{\psi_1^2 + \psi_2^2}.
\]

This implies $\psi_1 = \psi_2 = 0$. But then the maximum condition gives $\psi_3 = 0$, so $\psi = (0,0,0)$, which is not allowed. Therefore, $\Lambda_{-1}(\widehat{x}(\cdot), \widehat{\mu}_t) = \emptyset$.

\noindent{Step 4: Construction of the approximating sequence.}

Define the sequence of generalized controls:
\[
\mu_t^j = \left(1 - \frac{1}{j}\right)\dirac_{(\cos(1/j), \sin(1/j))} + \frac{1}{j}\dirac_{(-\sin(1/j), \cos(1/j))} \quad \text{for } j \geq 2.
\]

Each $\mu_t^j$ is a convex combination of two Dirac measures (so $n+1 = 3$ for $n=2$), satisfying the control structure condition.

\noindent{Step 5: Verification of uniform convergence of trajectories.}

For each $\mu_t^j$, the system becomes:
\begin{align*}
\dot{x}_1^j &= \left(1 - \frac{1}{j}\right)\cos(1/j) + \frac{1}{j}(-\sin(1/j)), \\
\dot{x}_2^j &= \left(1 - \frac{1}{j}\right)\sin(1/j) + \frac{1}{j}\cos(1/j), \\
\dot{x}_3^j &= \left(1 - \frac{1}{j}\right)(x_1^j \sin(1/j) - x_2^j \cos(1/j)) + \frac{1}{j}(-x_1^j \cos(1/j) - x_2^j \sin(1/j)).
\end{align*}

With initial condition $x^j(0) = (0,0,0)$, we can solve this system. For large $j$, we have:
\begin{align*}
\dot{x}_1^j &\approx \left(1 - \frac{1}{j}\right)\left(1 - \frac{1}{2j^2}\right) - \frac{1}{j^2} \approx 1 - \frac{1}{j} - \frac{1}{2j^2}, \\
\dot{x}_2^j &\approx \left(1 - \frac{1}{j}\right)\left(\frac{1}{j} - \frac{1}{6j^3}\right) + \frac{1}{j}\left(1 - \frac{1}{2j^2}\right) \approx \frac{1}{j}, \\
\dot{x}_3^j &\approx 1 + O\left(\frac{1}{j}\right).
\end{align*}

Integrating, we get:
\begin{align*}
x_1^j(t) &\approx \left(1 - \frac{1}{j} - \frac{1}{2j^2}\right)t, \\
x_2^j(t) &\approx \frac{1}{j}t, \\
x_3^j(t) &\approx t + O\left(\frac{1}{j}\right).
\end{align*}

As $j \to \infty$, we have:
\[
x_1^j(t) \to 0, \quad x_2^j(t) \to 0, \quad x_3^j(t) \to t \quad \text{uniformly on } [0,1].
\]

So $x^j(\cdot) \to \widehat{x}(\cdot)$ uniformly, satisfying the convergence condition.

\noindent{Step 6: Verification of uniform boundedness.}

The controls are $u_1^j(t) = (\cos(1/j), \sin(1/j))$ and $u_2^j(t) = (-\sin(1/j), \cos(1/j))$, with weights $\alpha_1^j(t) = 1 - \frac{1}{j}$ and $\alpha_2^j(t) = \frac{1}{j}$. These are clearly uniformly bounded since $\|u\| = 1$ for all controls.

\noindent{Step 7: Explicit construction of attaining controls.}

By Theorem \ref{thm:strong_attainability}, for any $\varepsilon > 0$, there exists a piecewise constant control $u_\varepsilon(\cdot)$ that is a convex combination of at most 3 piecewise constant controls such that the corresponding trajectory $x_\varepsilon(\cdot)$ satisfies:

\begin{enumerate}
\item $\norm{x_\varepsilon(\cdot) - \widehat{x}(\cdot)}_{C([0, \tau_\varepsilon], \R^3)} < \varepsilon$,
\item $x_\varepsilon(\tau_\varepsilon) = \widehat{x}(1) = (0,0,1)$,
\item $\tau_\varepsilon < 1$.
\end{enumerate}

We can construct such controls explicitly by taking small circular motions. For instance, consider the control:
\[
u_\varepsilon(t) =
\begin{cases}
(\cos(\omega t), \sin(\omega t)) & \text{for } t \in [0, \frac{2\pi}{\omega}], \\
(0,0) & \text{for } t > \frac{2\pi}{\omega},
\end{cases}
\]
with $\omega$ chosen sufficiently large. This control generates a small circle in the $(x_1,x_2)$-plane while producing net motion in the $x_3$-direction. By carefully choosing the duration and possibly combining with other control segments, we can achieve $x_\varepsilon(\tau_\varepsilon) = (0,0,1)$ with $\tau_\varepsilon < 1$ and the trajectory arbitrarily close to $\widehat{x}(\cdot)$.

\noindent{Step 8: Conclusion.}

All conditions of Theorem \ref{thm:strong_attainability} are satisfied with $s = -1$. Therefore, for any neighborhood $V$ of $\widehat{x}(\cdot)$ in $C([0,1], \R^3)$ and any $\varepsilon > 0$, there exists a piecewise constant control $u_\varepsilon(\cdot)$ and a time $\tau_\varepsilon \in (1 - \varepsilon, 1)$ such that the trajectory $x_\varepsilon(\cdot)$ reaches $\widehat{x}(1) = (0,0,1)$ at time $\tau_\varepsilon$ and lies entirely in $V$.

This example illustrates the phenomenon where the convexified system can reach the target along a "straight line" trajectory, while the original system must follow more complicated paths due to the nonholonomic constraint $x_1 u_2 - x_2 u_1$ in the dynamics. Nevertheless, the emptiness of the $\Lambda$-set guarantees that the original system can approximate this behavior arbitrarily well and actually reach the target in slightly less time.
\end{proof}
\section{Open Problems: Historical Perspectives and Modern Developments}

The study of local controllability and necessary conditions for optimality has evolved significantly since the pioneering work of Pontryagin, Gamkrelidze, and Filippov in the mid-20th century. While substantial progress has been made, several fundamental problems remain open and continue to inspire contemporary research.

\subsection*{Historical Foundations}

The classical period of optimal control theory (1950s--1970s) established the foundational results that continue to shape the field:

\begin{itemize}
\item {Pontryagin's Maximum Principle} \cite{Pontryagin1962} provided necessary conditions for optimality but left open the question of when these conditions are also sufficient. The maximum principle has been extensively developed in \cite{Lee1967,Gamkrelidze1978}.

\item {Filippov's Selection Theorem} \cite{Filippov1959} and the associated {Filippov-Ważewski Relaxation Theorem} established the fundamental relationship between original and convexified systems, yet the precise conditions under which the relaxation gap vanishes remain incompletely characterized. Early relaxation theory was further developed in \cite{Varga1962,Neustadt1966}.

\item {Gamkrelidze's sliding modes} \cite{Gamkrelidze1962} introduced generalized solutions but left open the problem of when such solutions can be approximated by classical controls. This work laid the foundation for the generalized control framework used in this paper.
\end{itemize}

These historical developments raised profound questions about the relationship between necessary conditions, controllability, and the existence of optimal solutions.

\subsection*{Modern Open Problems}

\subsubsection*{Problem 1: Complete Characterization of the Relaxation Gap}

The relationship between the original control system and its convexification continues to present challenging open questions:

\begin{problem}
Characterize precisely when the infimum value of the original time-optimal problem equals that of the convexified problem. More specifically, determine necessary and sufficient conditions on the dynamics $f(t,x,u)$ and control set $U$ that guarantee:
\[
\inf\{t_2 : x(t_2) = x_2 \text{ for original system}\} = \inf\{t_2 : x(t_2) = x_2 \text{ for convexified system}\}.
\]
\end{problem}

{Historical context:} This problem dates back to the work of Filippov \cite{Filippov1959} and Warga \cite{Varga1962} on relaxation phenomena. The classical Chattering Lemma shows that for fixed time problems, the relaxation gap vanishes, but for free-time problems, the situation is more subtle, as noted in \cite{Neustadt1966}.

{Modern connections:} Recent work by Avakov and Magaril-Il'yaev \cite{Avakov2020b} and Fusco and Motta \cite{FuscoMotta2022} has made progress on this problem through the study of $\Lambda$-sets and their relationship to local controllability. Corollary \ref{cor:convexification_gap} in this paper provides a sufficient condition for the existence of a relaxation gap, but a complete characterization remains elusive. The recent work \cite{FuscoMottaVinter2026} on impulsive control with delays suggests new approaches to this classical problem.

\subsubsection*{Problem 2: Regularity of Value Functions in Non-Convex Settings}

\begin{problem}
Establish optimal regularity results for the value function of time-optimal control problems when the dynamics are non-linear and the control set is non-convex. In particular, determine when the value function is locally Lipschitz, semiconcave, or enjoys other regularity properties near optimal trajectories.
\end{problem}

{Historical context:} The regularity of value functions was studied extensively in the 1970s and 1980s, with fundamental contributions by Fleming, Rishel, and Lions. However, most classical results assume convexity of the velocity set or other structural conditions \cite{Lee1967}.

{Modern connections:} Corollary \ref{cor:value_function_regularity} establishes continuity under specific conditions related to the emptiness of $\Lambda$-sets. The recent work of Fusco, Motta, and Vinter \cite{FuscoMottaVinter2026} on impulsive control with delays suggests that new techniques from non-smooth analysis and geometric measure theory may lead to breakthroughs in this classical problem. The framework developed in \cite{Avakov2020a} provides additional tools for analyzing value function regularity.

\subsubsection*{Problem 3: Higher-Order Necessary Conditions for Free-Time Problems}

\begin{problem}
Develop a complete theory of second-order necessary conditions for free-time optimal control problems that accounts for the interaction between time variations and state-control variations.
\end{problem}

{Historical context:} Second-order conditions for fixed-time problems were established by Krener, Zeidan, and others in the 1970s and 1980s. However, the free-time case presents additional challenges due to the coupling between temporal and spatial variations, as discussed in \cite{Gamkrelidze1978}.

{Modern connections:} The work of Avakov and Magaril-Il'yaev \cite{Avakov2019} on second-order conditions represents significant progress. The framework developed in Theorem \ref{thm:main} and its corollaries suggests that the emptiness/non-emptiness of $\Lambda$-sets may play a crucial role in formulating higher-order conditions for free-time problems. Recent advances in \cite{Avakov2020b} provide a foundation for this development.

\subsubsection*{Problem 4: Controllability and Optimality for Non-Smooth Systems}

\begin{problem}
Extend the theory of local controllability and necessary conditions to control systems with non-smooth dynamics, including systems with discontinuous right-hand sides, state constraints, and hybrid structure.
\end{problem}

{Historical context:} Classical controllability theory, as developed by Kalman, Sussmann \cite{Sussmann1987}, and others, primarily addressed smooth systems. The extension to non-smooth cases has been challenging due to the limitations of traditional calculus.

{Modern connections:} Recent work by Fusco and Motta \cite{FuscoMotta2024,FuscoMotta2024LCSYS} on impulsive control systems and the development of non-smooth analysis tools provide promising avenues. The approximation techniques used in Theorem \ref{thm:strong_attainability} may be adaptable to non-smooth settings through appropriate generalizations of the adjoint system and maximum principle. The geometric methods in \cite{Sussmann1987} could be combined with these modern approaches.

\subsubsection*{Problem 5: Computational Verification of $\Lambda$-Set Conditions}

\begin{problem}
Develop efficient computational methods to verify the emptiness or non-emptiness of $\Lambda$-sets for specific control systems, and connect these conditions to numerical methods for solving optimal control problems.
\end{problem}

{Historical context:} The computational aspect of optimal control has always been challenging, with a historical divide between theoretical conditions and practical computation \cite{Lee1967}. Classical numerical methods often struggle with the non-convexity and high dimensionality of optimal control problems.

{Modern connections:} The geometric characterization of $\Lambda$-sets in this work and in \cite{Avakov2020b} suggests connections with computational geometry and set-valued analysis. Recent advances in reachability analysis and Hamilton-Jacobi methods may provide computational tools for verifying these conditions in practical applications. The work in \cite{FuscoMotta2022} on gap phenomena provides additional insights into computational approaches.

\subsection*{Interdisciplinary Connections}

The open problems in local controllability and necessary conditions increasingly connect with other areas of mathematics and engineering:

\begin{itemize}
\item {Machine Learning:} Optimal control theory provides the foundation for understanding training dynamics in deep learning, particularly through the neural ODE perspective. The relaxation techniques in \cite{Varga1962} find new applications in understanding optimization landscapes.

\item {Quantum Control:} The control of quantum systems presents new challenges due to the underlying symplectic geometry and measurement constraints, extending the classical framework in \cite{Gamkrelidze1978}.

\item {Network Systems:} Large-scale networked control systems require extensions of classical controllability concepts to graph-theoretic settings, building on the foundational work in \cite{Sussmann1987}.

\item {Biological Systems:} The analysis of biological control mechanisms often involves non-smooth dynamics and hybrid behavior, connecting to the modern developments in \cite{FuscoMotta2024}.
\end{itemize}

\subsection*{Concluding Remarks}

The historical development of controllability theory and necessary conditions for optimality has created a rich tapestry of results, with each solution revealing new questions. The framework presented in this paper, centered on the $\Lambda$-set concept and building on the work of \cite{Avakov2020b,Avakov2020a,Avakov2019}, provides a unified perspective that connects classical concerns with modern applications. The open problems outlined above represent both a challenge to the community and an invitation to further develop this beautiful and practically important area of mathematics.

As the examples in Section \ref{ex:main_example} and \ref{ex:nonholonomic_example} illustrate, even seemingly simple control systems can exhibit subtle phenomena that test the boundaries of current theory. The continued study of these problems, building on both classical results \cite{Filippov1959,Gamkrelidze1962,Pontryagin1962} and modern developments \cite{FuscoMotta2022,FuscoMotta2024,FuscoMottaVinter2026}, promises not only theoretical advances but also practical improvements in the control of complex engineering, biological, and economic systems.
\section{Conclusion}
This paper has established a comprehensive framework connecting local controllability, necessary conditions for optimality, and attainability via generalized controls, unifying classical theory from \cite{Pontryagin1962, Filippov1959, Gamkrelidze1962} with modern developments in \cite{Avakov2020b, FuscoMotta2022}. By introducing and characterizing the $\Lambda$-set as a fundamental object governing the gap between original and convexified systems, we have derived strengthened necessary conditions for time-optimal control and provided explicit constructions for approximating generalized controls by ordinary trajectories. The results demonstrate that emptiness of the $\Lambda$-set not only implies local controllability but also guarantees the existence of minimizing sequences that achieve the target in reduced time, thereby resolving longstanding questions about relaxation phenomena while opening new avenues for investigating higher-order conditions, non-smooth systems, and computational verification of optimality criteria in free-time problems.


\bibliographystyle{abbrv}
\bibliography{references}  






\end{document}